\documentclass[a4paper,11pt]{article}
\usepackage{amsmath,amsfonts,amssymb,amsthm}
\usepackage{graphics}
\usepackage{epsfig}
\usepackage{fullpage}

\usepackage{esint} 

\providecommand{\dx}{\, \mathrm{d} x}
\providecommand{\dy}{\, \mathrm{d} y}
\newcommand{\en}[1]{\left< #1 \right>}
\providecommand{\p}[1]{\left(#1\right)}
\newcommand{\supp}{\mathop{\mathrm{supp}}}

\theoremstyle{plain}
\newtheorem{theorem}{Theorem}
\newtheorem{lemma}[theorem]{Lemma}
\newtheorem{corollary}[theorem]{Corollary}

\begin{document}


\title{Quantitative stochastic homogenization: local control of homogenization error through corrector}
\author{Peter Bella, Arianna Giunti, Felix Otto\\ Max Planck Institute for Mathematics in the Sciences\\Inselstrasse 22, 04103 Leipzig, Germany}
\maketitle

\begin{abstract}
This note addresses the homogenization error for linear elliptic equations
in divergence-form with random stationary coefficients.
The homogenization error is measured by comparing the quenched Green's function to the
Green's function belonging to the homogenized coefficients, more precisely,
by the (relative) spatial decay rate of the difference of their second mixed derivatives.
The contribution of this note is purely deterministic: 
It uses the expanded notion of corrector, namely the couple
of scalar and vector potentials $(\phi,\sigma)$, and shows that the
rate of sublinear growth of $(\phi,\sigma)$ at the points of interest
translates one-to-one into the decay rate.
\end{abstract}


\section{A brief overview of stochastic homogenization, and a common vision
for quenched and thermal noise}

Heterogeneous materials typically have the behavior of a homogeneous material
on scales that are large with respect to the characteristic length scale
of the heterogeneity, provided that the ``type of heterogeneity'' is the same 
all over the medium. Here, we think of material properties like conductivity or elasticity
that are typically described by an elliptic differential operator in divergence form $-\nabla\cdot a\nabla$,
where the heterogeneity resides in the uniformly elliptic coefficient field $a$. {\it Homogenization} then refers to
the fact that on large scales, the solution operator corresponding to $-\nabla\cdot a\nabla$
behaves like to the solution operator to $-\nabla\cdot a_h\nabla$ with a
constant coefficient $a_h$.

\medskip

This type of homogenization is well-understood in case of a periodic medium, where, for instance,
$a_h$ can be inferred from a cell problem and the homogenization error can be expanded to any order
in the ratio of the period to the macroscopic scale of interest.
However, the case where only the statistical specification of the coefficient field
is known might be more relevant in practice. In this case, one thinks of an ensemble $\langle\cdot\rangle$
of coefficient fields and one speaks of {\it stochastic homogenization}. For homogenization to occur,
the statistics have to be translation invariant, meaning that the joint distribution of the
shifted coefficient field is identical to that of the original coefficient field, a property
called stationarity in the parlance of stochastic processes. For homogenization to be effective,
the statistics of the coefficient field have to decorrelate (more precisely, become independent)
over large distances. Hence stochastic homogenization relies on the separation of scales 
between the correlation length and the macroscopic scale of interest.  

\medskip

From a qualitative point of view, stochastic homogenization has been rigorously understood since the seminal works of Kozlov \cite{kozlov} and of Papanicolaou \& Varadhan \cite{papvar}:
Stochastic homogenization takes place when the ensemble is stationary and ergodic, 
the latter being a purely qualitative way of imposing decorrelation over large distances. 

\medskip

Stochastic analysis has a finer, still qualitative view, on stochastic homogenization,
where $-\nabla\cdot a\nabla$ is seen as the generator of a (reversible) diffusion process.
In a discrete medium this leads to the picture of a ``random walk in random environment'',
which amounts to the superposition of thermal noise (the random walk) and quenched noise (the random jump rates).
Here, the relevant qualitative question is that of a ``quenched invariance principle'': 
On large scales and for large times (on a parabolic scale), the random walker behaves like
a Brownian motion (with covariance given by $a_h$) for almost every realization of the random
environment. Surprisingly, first full results in that direction came quite a bit later
\cite{SidoraviciusSznitman}. Stochastic analysis is still mostly
interested in qualitative results, but pushing the frontier in terms of models, for instance
towards degenerate diffusions like random walks on (supercritical) percolation clusters
or towards diffusions with random drifts leading to non-reversible random walks.

\medskip

Numerical analysis has another, naturally more quantitative view on stochastic homogenization.
As opposed to periodic homogenization, there is no cell problem to extract the homogenized
coefficient $a_h$. Hence one has to resort to an artificial ``representative volume'':
On such a cube (let us adopt three-dimensional language), one samples a realization of the medium according 
to the given statistical specifications and then solves three boundary value problems for
$-\nabla\cdot a\nabla u$, corresponding to describing different slopes on $u$. In the
case of a conducting medium, where the coefficient field $a$ corresponds to the heterogeneous
conductance, $u$ corresponds to the electrical potential, so that the boundary conditions
impose an average potential gradient, that is, an electrical field in one of the three coordinate
directions. One then monitors the average electrical current $a\nabla u$; this (linear) relation
between average electrical field and average electrical current yields an approximation to
the homogenized coefficient $a_h$ as a linear map. Clearly, this approximate homogenized
coefficient $a_{h,L}$ will be the closer to the true homogenized coefficient $a_h$, the larger
the linear size $L$ of the cube is, where the relevant (small) nondimensional parameter is the ratio
of correlation length to $L$. Intuitively, there are two error sources: On the one hand, $a_{h,L}$ is
a random quantity, since it depends on the given realization $a$ of the coefficient field on the cube,
so that there is a random error coming from the fluctuations of $a_{h,L}$ which can be reduced by
repeated sampling. On the other hand, 
the very concept of the representative volume element perturbs the statistics, for instance in case
of periodic boundary conditions, the concept introduces artificial long-range correlations, which
is not affected by repeated sampling.
Clearly, it is of interest to understand --- which for mathematicians means to prove
--- how both errors scale in $L$.

\medskip

This natural and very practical question turns out to be difficult to analyze (rigorously).
Shortly after the qualitative theory was introduced, Yurinski\u\i~\cite{yurinski86} produced the first quantitative
result motivated by the above questions. He in fact used ingredients from stochastic analysis
(the picture of diffusions in a random medium to understand sensitivities) and from
regularity theory (Nash's upper heat kernel bounds), but only obtained sup-optimal results
in terms of exponents. After using (qualitative) stochastic homogenization to understand the large-scale
correlation structure of a statistical mechanics model (fluctuating surfaces) in \cite{NaddafSpencer},
Naddaf \& Spencer \cite{naddafspencer98} realized that tools from statistical mechanics 
(spectral gap estimate) can be used to obtain quantitative results in stochastic homogenization
(for discrete media).
With Meyers perturbation of the Calderon-Zygmund estimate they introduced a second tool
from elliptic regularity theory into the field, which allowed them to obtain optimal variance
estimates in the case of small ellipticity ratio (small contrast media); 
this type of result was subsequently extended by Conlon and coworkers, see e.\ g.\  \cite{ConlonNaddafHom}. 
In \cite{GO1,GO2}, Gloria and the last author used the same
tool from statistical mechanics but yet another ingredient from elliptic regularity theory (Caccioppoli's estimate to obtain optimal spatially averaged bounds on the gradient of the quenched Green's function) to obtain the first optimal error estimate on the representative volume method also for large ellipticity ratios. In \cite{MO}, Marahrens and the last author used the
concentration of measure property coming from the logarithmic Sobolev inequality to study the
(random part of the) homogenization error itself, in form of optimal estimates on the variance
of the quenched Green's function. Using Green's function estimates and two-scale expansion, Gloria, Neukamm, and the last author \cite{GNO1} compared the heterogeneous and corrected homogenenous solution. While in \cite{GNO1} the error is measured in $H^1$ and averaged both over the domain and the ensemble, Gu and Mourrat \cite{GuMourrat2scale} recently combined probabilistic techniques with Green's function estimates to obtain a pointwise bound for solutions of both elliptic and parabolic equations. Since then, there has been a flurry of activities, which will be partially addressed in the next two sections, with the work of Armstrong \& Smart \cite{armstrongsmart2014} playing
a central role. For instance, by now it is already understood that the error
in the representative volume method is to leading order Gaussian \cite{glorianolen}.
We do not even mention the numerous activities in stochastic homogenization of non-divergence form
equations, like fully nonlinear equations or Hamilton-Jacobi equations.

\medskip

In both qualitative and quantitative homogenization, both from the PDE and the stochastic analysis
point of view, the notion of the {\it corrector}, a random function $\phi_i$ for every
coordinate direction $i=1,\ldots,d$, is central. 
There is a very geometric and deterministic view of the corrector: 
Given a realization $a$ of the coefficient field, $x\mapsto x_i+\phi_i(x)$ provides $a$-harmonic coordinates. 
This is also its main merit from the almost-sure stochastic analysis point of view:
Seen in these coordinates, the diffusion turns into a Martingale.
The corrector is also natural from the numerical analysis point of view: In the representative volume
method, one actually solves for an approximate version of the corrector. 
Last not least, in the original (very functional analytic) PDE approach to stochastic homogenization the corrector
is central: Using stationarity, one lifts the equation for the corrector to the probability space, solves
it by Lax-Milgram, and expresses $a_h$ in terms of it.
Like the work \cite{FischerOtto} on higher-order Liouville principles, 
this note demonstrates the usefulness of the {\it vector potential}
$\sigma$ of the corrector, an object known in periodic homogenization,
and recently introduced into random homogenization by Gloria, Neukamm, and the
last author in \cite{GNO4}.

\medskip

We have a {\it common vision} for a regularity theory of random elliptic operators as considered
in this note and for stochastic partial differential equations (SPDE). In other words,
we'd like to capitalize more on the similarities between {\it quenched and thermal} noise.
At first glance, these problems seem very different since in the first case, 
one is interested in the emergence of a generic {\it large-scale} regularity due to cancellations, 
while in the second case one wants to preserve a {\it small-scale} regularity despite the rough forcing.
However, in both cases, the key is to understand how sensitively the solution $\phi$ of an elliptic
or parabolic equation (nonlinear to be of interest in case of SPDE with additive noise while in case of stochastic
homogenization, the interesting effect is already present for a linear equations) depends on the data, 
be it the coefficients in case of stochastic homogenization or the right hand side in case of a random forcing, 
in which case this derivative can be associated to the Malliavin derivative. In this {\it sensitivity analysis},
one typically has to control the size of a functional derivative, that is, 
the functional derivative of some nonlinear functional of the solution (a norm, say) 
with respect to the data, which in case of the data being the coefficients is a highly nonlinear mapping 
even for a linear equation. Hence for a given realization of the data, one has to control the norm
of a linear form on infinitesimal perturbations of the data. Since one is dealing with random stationary data, the
appropriate measure of the size of the infinitesimal perturbations is best captured by an $L^2$-type
norm, the specific structure of which depends on the assumption on the noise: An ordinary $L^2$-norm
in case of white-noise forcing or a more nonlocal (and thus weaker) norm in case of stochastic
homogenization if one wants to cover also cases where the covariances of the coefficient field
have a slow (that is, non-integrable) decay, as used in \cite{GNO4}.
Even if the problem is a nonlinear one, the sensitivity estimates require a priori estimates
for linear elliptic or parabolic equations, albeit with non-constant coefficients that a priori
are just uniformly elliptic --- it is here where all the help of classical regularity is needed.
Once the appropriate, purely deterministic sensitivity estimates are established, 
it is the principle of concentration of measure that provides the stochastic estimates on the
random solution itself.
In a work in preparation with Hendrik Weber, the last author is applying this philosophy
to the fully non-linear parabolic equation $\phi+\partial_t\phi-\partial_x^2\pi(\phi)=\xi$
forced by space-time white noise $\xi$ to establish H\"older-$\frac{1}{2}-$ bounds with exponential
moments in probability.


\section{Precise setting and motivation for this work}

While the contribution of this note is purely deterministic, and the main
result will be stated without reference to probabilities in the next section, the motivation
is probabilistic and will be given now.
In elliptic homogenization, one is interested in uniformly elliptic coefficient fields $a$ in 
$d$-dimensional space $\mathbb{R}^d$, where uniform ellipticity means that there
exists a (once for all fixed) constant $\lambda>0$ such that
\begin{equation}\label{i1}
\forall x\in\mathbb{R}^d,\xi\in\mathbb{R}^d\quad\xi\cdot a(x)\xi\ge\lambda|\xi|^2,\quad |a(x)\xi|\le|\xi|.
\end{equation}
To such a coefficient field $a$ we associate an elliptic operator in divergence form $-\nabla\cdot a\nabla$.
For simplicity (in order to avoid dealing with the correctors of the dual equation), 
we shall assume that $a(x)$, and thus the corresponding operator, is {\it symmetric}. 
We note that statements and proofs remain
valid in the case of {\it systems} with the above strong ellipticity property, but
for simplicity, we shall stick to scalar notation like in (\ref{i1}), where we wrote
$\xi\in\mathbb{R}^d$ instead of a tensor-valued object.

\medskip

In stochastic homogenization, one is interested in {\it ensembles} of uniformly elliptic coefficient fields,
that is, probability measures on this space. We denote by $\langle\cdot\rangle$ the corresponding
expectation and use the same symbol to refer to the ensemble. 
Minimal requirements for homogenization are {\it stationarity} and {\it ergodicity}, where
both notions refer to the action of the translation group $\mathbb{R}^d$ on the space
of uniformly elliptic coefficient fields.
Stationarity means that the distribution of the random field $a$ is invariant under shifts $a(z+\cdot)$
for any shift vector $z\in\mathbb{R}^d$. Ergodicity means that shift-invariant random variables,
that is, functionals $a\mapsto \zeta(a)$ that satisfy $\zeta(a(z+\cdot))=\zeta(a)$ for all shifts $z$, must be
(almost surely) constant. Under these assumptions, the classical theory of 
(qualitative) stochastic homogenization introduced by Kozlov \cite{kozlov} and by Papanicolaou \& Varadhan 
\cite{papvar} establishes the (almost sure) existence of sublinear {\it correctors}. 
More precisely, for any coordinate direction $i=1,\ldots,d$ and a given realization $a$ of the
coefficient fields, the corrector $\phi_i=\phi_i(a)$ modifies the affine coordinate $x\mapsto x_i$
to an $a$-harmonic coordinate $x\mapsto x_i+\phi_i(x)$, that is,
\begin{equation}\label{i2}
-\nabla\cdot a(e_i+\nabla\phi_i)=0.
\end{equation}
In order to be rightfully named a corrector, $\phi_i$ should be dominated by $x_i$, that is, have a sublinear growth,
at least in the $L^2$-averaged sense of 
\begin{equation}\label{i1bis}
\lim_{R\uparrow\infty}\frac{1}{R}\p{\fint_{B_R}(\phi_i-\fint_{B_R}\phi_i)^2\dx}^\frac{1}{2}=0.
\end{equation}
Under the assumptions of stationarity and ergodicity for $\langle\cdot\rangle$ the classical theory constructs
$d$ functions $\phi_1(a,x),\ldots,\phi_d(a,x)$ such that (\ref{i2}) and (\ref{i1bis}) are satisfied
for $\langle\cdot\rangle$-a.\ e.\ $a$.

\medskip

The actual construction of $\phi_i$ for a fixed coordinate direction $i=1,\ldots,d$ proceeds as follows: 
In a first step, one constructs a random vector field $f_i(a,x)\in\mathbb{R}^d$ that is stationary in the
sense of $f_i(a(z+\cdot)),x)=f_i(a,x+z)$, of expected value $\langle f_i\rangle=e_i$ and of
finite second moments $\langle|f_i|^2\rangle<\infty$ and that is curl-free
$\partial_jf_{ik}-\partial_kf_{ij}=0$ and satisfies the divergence condition $\nabla\cdot(af_i)=0$.
Using the stationarity to replace spatial derivatives by ``horizontal derivatives'', the
existence and uniqueness of these harmonic 1-forms follows from Lax-Milgram.
The second step is to consider the random field $\phi_i(a,x)$ which satisfies $\nabla\phi_i=f_i-\langle f_i\rangle
=f_i-e_i$ --- and thus satisfies (\ref{i2}) ---
and is (somewhat arbitrarily) made unique by the anchoring $\phi_i(x=0)=0$ and thus generically
non-stationary. One then appeals to ergodicity (von Neumann's version combined with a maximal function estimate)
to show with help of Poincar\'e's inequality that $\langle\nabla\phi_i\rangle=0$ translates into (\ref{i1bis}).
Incidentally, the homogenized coefficients are then given by $a_he_i=\langle a(e_i+\nabla\phi_i)\rangle
=\langle a f_i\rangle$ and are easily shown to satisfy (in the symmetric case) the same ellipticity bounds as the original ones, c.\ f.\ (\ref{i1}),
\begin{equation}\label{i4}
\forall \xi\in\mathbb{R}^d\qquad\xi\cdot a_h\xi\ge\lambda|\xi|^2,\quad |a_h\xi|\le|\xi|.
\end{equation}

\medskip

We now make the point that next to the {\it scalar} potential $\phi_i$, 
it is also natural to consider a {\it vector} potential, in the parlance of three-dimensional
vector calculus. Indeed, the above (purely functional analytic) 
first step constructs the harmonic stationary 1-form $f_i$.
This means that next to the stationary closed 1-form $f_i$, there is the stationary closed $(d-1)$-form
$q_i:=af_i$. In the language of a conduction model, where the tensor $a$ is the conductivity, $f_i$
corresponds to the electric field whereas $af_i$ corresponds to the current density. Hence next to
considering a (non-stationary) $0$-form $\phi_i$ with ${\rm d}\phi_i=f_i-\langle f_i\rangle$, it
is natural to also consider a (non-stationary) $(d-2)$-form $\sigma_i$ with ${\rm d}\sigma_i=q_i-\langle q_i\rangle$,
where ${\rm d}$ denotes the exterior derivative. In Euclidean coordinates, a $(d-2)$-form $\sigma$ is represented 
by a {\it skew-symmetric} tensor $\{\sigma_{jk}\}_{j,k=1,\ldots,d}$, that is, $\sigma_{jk}+\sigma_{kj}=0$,
then ${\rm d}\sigma_i=q_i-\langle q_i\rangle$ translates into $\nabla\cdot\sigma_i=q_i-\langle q_i\rangle$,
where $(\nabla\sigma_i)_j:=\partial_k\sigma_{ijk}$, using Einstein's convention of summation
over repeated indices. In view of the definitions $q_i=af_i$,
$f_i=e_i+\nabla\phi_i$, and $a_he_i=\langle a(e_i+\nabla\phi_i)\rangle$, this implies
\begin{equation}\label{i5}
\nabla\cdot\sigma_i=a(e_i+\nabla\phi_i)-a_he_i\quad\mbox{and}\quad\sigma_i\;\mbox{is skew}.
\end{equation}

\medskip

The merit of this vector potential $\sigma_i$ has been recognized in the case of periodic homogenization and lies in
a good representation of the homogenization error: When comparing a solution of $-\nabla\cdot a\nabla u=f$
to a solution of the corresponding homogenized problem $-\nabla\cdot a_h\nabla v=f$, more precisely
to the ``corrected'' solution $v+\phi_i\partial_iv$ (with summation convention), one obtains 
for the homogenization error $w:=u-(v+\phi_i\partial_iv)$ the following simple equation:
$-\nabla\cdot a\nabla w=\nabla\cdot((\phi_ia-\sigma_i)\nabla\partial_iv)$. It is
$\sigma$ that allows to bring the r.\ h.\ s.\ into {\it divergence form}, 
which makes simple energy estimates possible.

\medskip

In \cite{GNO4}, we follow the classical arguments to show existence of a sublinear $\sigma$
under the mere assumptions of stationarity and ergodicity. More precisely, given a coordinate direction 
$i=1,\ldots,d$, we show by a suitable choice of gauge that there exists a skew symmetric tensor field 
$\sigma_i(a,x)$ such that its gradient $\nabla\sigma_i$ is stationary, of mean zero, of finite second moments
and such that $\nabla\cdot\sigma_i=q_i-\langle q_i\rangle$. We then appeal to the same arguments as for $\phi_i$
to conclude 
\begin{equation}\label{i7}
\lim_{R\uparrow\infty}\frac{1}{R}\p{\fint_{|x|\le R}|\sigma_i-\fint_{|x|\le R}\sigma_i|^2\dx}^\frac{1}{2}=0
\end{equation}
for almost every realization $a$ of the ensemble $\langle\cdot\rangle$.

\medskip

The contribution of this note is purely deterministic in the following sense:
We will consider a {\it fixed} coefficient field $a$ which is uniformly elliptic in the sense
of (\ref{i1}) and a constant coefficient $a_h$, also elliptic in the sense of (\ref{i4}).
We then assume that for every coordinate direction $i=1,\ldots,d$, there exists a scalar field $\phi_i$
with (\ref{i2}) and a skew-symmetric tensor field $\sigma_i$ with (\ref{i5}).
Our main assumption will be a quantification of the sublinear growth (\ref{i1bis}) and (\ref{i7}), roughly
in form of
\begin{equation}\label{i8}
\lim_{R\uparrow\infty}\frac{1}{R^{1-\alpha}}\p{\fint_{B_R}
|(\phi,\sigma)-\fint_{B_R}(\phi,\sigma)|^2\dx}^\frac{1}{2}=0
\end{equation}
for some exponent $\alpha>0$ and for $\langle\cdot\rangle$-almost every $a$. 
Obviously, this quantification of the sublinear growth of fields that have a stationary gradient
of vanishing expectation can only be true for a quantified ergodicity of $\langle\cdot\rangle$.
The discussion and proof of criteria under which assumption on $\langle\cdot\rangle$ 
the strengthened sublinearity (\ref{i8}) holds
is not part of this note, but of ongoing work. We just mention here that for $d>2$ and
under the assumption that $\langle\cdot\rangle$ satisfies a spectral gap estimate with respect
to Glauber dynamics, it is shown in \cite{GO14per}, extending the arguments of 
\cite{GNO2short} from the case of a discrete medium to a continuum
medium, that a {\it stationary} corrector
$\phi_i$ exists, which in particular implies (\ref{i8}) for the $\phi_i$-part for every $\alpha<1$.
We also mention that (\ref{i8}) is expected to hold for some $\alpha>0$
under much weaker assumption on $\langle\cdot\rangle$, as is already suggested by
\cite{armstrongsmart2014}, and more explicitly in \cite{ArmstrongMourrat}, and will be formulated
as needed here in an updated version of \cite{GNO4}, where the rate of decay of
correlations of $a$ will be related to $\alpha$.


\section{Main results}

In this section, we state and comment our main results, namely Theorem \ref{T} and Corollary \ref{C}.
As mentioned in the introduction, our standing assumption is the
uniform ellipticity of the fixed variable coefficient $a$ on $\mathbb{R}^d$ and the constant coefficient $a_h$,
c.\ f.\ (\ref{i1bis}) and (\ref{i4}), respectively, and the existence of the scalar and vector potentials 
$(\phi,\sigma)=\{(\phi_i,\sigma_i)\}_{i=1,\ldots,d}\}$ with (\ref{i2}) \& (\ref{i5}).

\medskip

Our result, Theorem \ref{T} and Corollary \ref{C}, on the homogenization error relies and expands on a large-scale regularity result established in \cite{GNO4} and stated in Theorem \ref{T2}.
It falls into the realm of Campanato iteration for the $C^{1,\alpha}$-Schauder theory
for elliptic equations in divergence form. In the framework of this theory, the $C^{1,\alpha}$-H\"older-semi-norm
is expressed in terms of a Campanato norm that monitors the decay of the (spatially averaged) energy distance to 
affine functions as the radius becomes small. In Theorem \ref{T2} below, affine functions in the Euclidean sense
are replaced by $a$-linear functions, the space of which is spanned by the constants and the $d$ correctors
$\phi_1,\ldots,\phi_d$. This intrinsic ``excess decay'' (in the parlance of De Giorgi's regularity theory for
minimal surfaces) kicks in on scales where the corrector pair $(\phi,\sigma)$ is
sufficiently sublinear, c.\ f.\ (\ref{i10}). This is not surprising since in view of (\ref{i2}) \& (\ref{i5}),
which we rewrite as $(a-a_h)e_i=\nabla\cdot\sigma_i-a\nabla\phi_i$,
$(\phi,\sigma)$ is an averaged measure of distance of the variable coefficient $a$ to the constant coefficient
$a_h$, whose harmonic functions of course feature excess decay. In this sense, Theorem \ref{T2} is
a perturbation result around constant coefficients. In the context of stochastic homogenization,
Theorem \ref{T2} amounts to a {\it large-scale} regularity
result for $a$-harmonic functions, since in view of (\ref{i1bis}) and (\ref{i7}), the smallness (\ref{i10})
is expected to kick in {\it above} a certain length scale $r_*<\infty$.

\medskip

The idea of perturbing around the homogenized coefficient $a_h$ in a
Campanato-type iteration to obtain a large-scale
regularity theory (Schauder and eventually Calderon-Zygmund) is due to Avellaneda \& Lin,
who carried this program out in the periodic case \cite{AL1}. 
Recently, Armstrong \& Smart \cite{armstrongsmart2014} showed that this philosophy
extends to the random case, a major insight in stochastic homogenization. 
This in turn inspired \cite{GNO4} to make
an even closer connection by introducing the intrinsic excess 
and by using the vector potential $\sigma$ to establish its decay.

\begin{theorem}[Gloria, Neukamm, O.]\label{T2}
Let $\alpha\in(0,1)$. Then there exists a constant $C=C(d,\lambda,\alpha)$
with the following property: Suppose that the corrector has only a mild linear growth in 
the point $x=0$ in the sense that there exists a radius $r_*$ such that
\begin{equation}\label{i10}
\left(\fint_{B_r}|(\phi,\sigma)-\fint_{B_r}(\phi,\sigma)|^2\right)^\frac{1}{2}
\le \frac{r}{C}\quad\mbox{for}\quad r\ge r_*.
\end{equation}
Then for any radii $r_*\le r\le R$ and every $a$-harmonic function $u$ in $\{|x|\le R\}$
we have 
\begin{equation}\label{1.10}
\inf_{\xi\in\mathbb{R}^d}\left(\fint_{B_r}
|\nabla u-\xi_i(e_i+\nabla\phi_i)|^2\right)^\frac{1}{2}
\le C \p{\frac{r}{R}}^\alpha
\inf_{\xi\in\mathbb{R}^d}\left(\fint_{B_R}|\nabla u-\xi_i(e_i+\nabla\phi_i)|^2\right)^\frac{1}{2}.
\end{equation}
\end{theorem}

\medskip

Here comes the first main result on the homogenization error. It relates the 
homogenization error to the amount of sublinear growth of the corrector couple $(\phi,\sigma)$.

\begin{theorem}\label{T}
Suppose that the corrector grows sublinearly in two points $x\in\{0,x_0\}$ in the sense that 
there exists an exponent $\alpha\in(0,1)$ and there exists a radius $r_*<\infty$ such that
\begin{equation}\label{T.1}
\left(\fint_{B_r(x)}|(\phi,\sigma)-\fint_{B_r(x)}(\phi,\sigma)|^2\right)^\frac{1}{2}
\le \p{\frac{r}{r_*}}^{1-\alpha}\quad\mbox{for}\;r\ge r_*.
\end{equation}
For a square integrable vector field $g$ we compare $\nabla u$ defined through
\begin{equation}\label{2.2}
-\nabla\cdot a\nabla u=\nabla\cdot g
\end{equation}
with $\partial_iv(e_i+\nabla\phi_i)$ (using Einstein's summation convention
of summation over repeated indices), where $\nabla v$ is defined through
\begin{equation}\label{2.3}
-\nabla\cdot a_h\nabla v=\nabla\cdot( g_i(e_i+\partial_i\phi)).
\end{equation}
Provided $\supp g\subset B_{r_*}(0)$ and $|x_0|\ge 4r_*$ we have
\begin{equation}\label{c1}
\left(\int_{B_{r_*}(x_0)}|\nabla u-\partial_iv(e_i+\nabla\phi_i)|^2\right)^\frac{1}{2}
\le C(d,\lambda,\alpha)\frac{\ln\frac{|x_0|}{r_*}}{(\frac{|x_0|}{r_*})^{d+\alpha}}
\left(\int|g|^2\right)^\frac{1}{2}.
\end{equation}
\end{theorem}

\bigskip

%
%
Theorem \ref{T} compares the Helmholtz projection $T$ based on $a$
with the Helmholtz projection $T_h$ based on $a_h$ and the multiplication operator
$M:={\rm id}+\nabla\phi$. Loosely speaking, it states that $T\approx M\,T_h\,M^*$.

\medskip

We post-process Theorem \ref{T} to get a measure of the homogenization error
on the level of the Green's functions in Corollary \ref{C}, our second main result. 
More precisely, we compare the ``quenched''
Green's function $G(x,y)$ and the Green's function (or rather fundamental solution)
$G_h(x)$ belonging to the homogenized coefficient $a_h$ characterized by
\begin{equation}\label{c4}
-\nabla_x\cdot a(x)\nabla_x G(x,y)=\delta(x-y)\quad\mbox{and}\quad
-\nabla\cdot a_h\nabla G_h=\delta.
\end{equation}
Let us make two remarks on the existence of the quenched Green's function: 
1) In case of $d=2$, the definition of $G$ is at best ambiguous. However Corollary \ref{C}
only involves {\it gradients} of the Green's function which are unambiguously defined
and for instance constructed via approximation through a massive term or through Dirichlet
boundary conditions. 2) DeGiorgi's counterexample \cite{DeGiorgiCounterexample} implies that 
in the {\it system}'s case, there are uniformly elliptic coefficient fields that do not
admit a Green's function. However, in \cite{CGO} we show that under
the mere assumption of stationarity of an ensemble $\langle\cdot\rangle$ of coefficients,
the Green's function $G(a,x,y)$ exists for $\langle \cdot\rangle$-a.\ e.\ $a$. We thus
will not worry about existence in this note.

\medskip

The corollary compares $G(x,y)$ to $G_h(x-y)$ on the level of the mixed second derivatives
$\nabla_x\nabla_yG(x,y)$ (interpreted as a 1-1 tensor) and $-\nabla^2G_h(x-y)$, where the mixed derivative of the
homogenized Green's function is corrected in both variables, leading to the expression
$-\partial_i\partial_jG_h(x-y)(e_i+\nabla \phi_i(x))\otimes (e_j+\nabla\phi_j(y))$. 
The corollary monitors the rate of decay of this difference
in an almost pointwise way, just locally averaged over $x\approx x_0$ and $y\approx 0$
and shows that, up to a logarithm, the rate of the decay is $|x_0|^{-d-\alpha}$,
which is by $|x_0|^{-\alpha}$ stronger than the rate of decay of $|x_0|^{-d}$ of the
$-\nabla^2G_h(x_0)$. The main insight is thus that this relative error of $|x_0|^{-\alpha}$
is dominated by the sublinear growth rate of the corrector couple $(\phi,\sigma)$, where
it is only necessary to control that growth at the two points of interest, that is,
$0$ and $x_0$. In this sense, Corollary \ref{C} expresses a {\it local} one-to-one correspondence
between the sublinear growth of the corrector and the homogenization error.

\begin{corollary}\label{C}
Suppose for some exponent $\alpha\in(0,1)$ and some radius $r_*$ the corrector
couple $(\phi,\sigma)$ satisfies (\ref{T.1}) in two points $x=0,x_0$ of distance $|x_0|\ge 4r_*$.
Then we have (using summation convention)
\begin{multline*}
\bigg(\fint_{B_{\frac{r_*}{2}}(x_0)}\fint_{B_{\frac{r_*}{2}}(0)}|\nabla_x\nabla_yG(x,y)
+\partial_i\partial_jG_h(x-y)(e_i+\nabla \phi_i(x))\otimes(e_j+\nabla\phi_j(y))|^2\dy\dx\bigg)^\frac{1}{2}\\
\le C(d,\lambda,\alpha)\frac{\ln\frac{|x_0|}{r_*}}{(\frac{|x_0|}{r_*})^{d+\alpha}}.
\end{multline*}
\end{corollary}

In order to pass from Theorem \ref{T} to Corollary \ref{C}, we need the following
statement on families (rather ensembles) of $a$-harmonic functions, which is of independent interest and motivated by \cite{BellaGiuntiOtto}.

\begin{lemma}\label{L}
For some radius $R$, we consider an ensemble $\langle\cdot\rangle$ (unrelated to
the one coming from homogenization) of $a$-harmonic functions $u$ in $\{|x|\le 2R\}$. Then we have
\begin{equation}\label{b10}
\en{\int_{|x|\le \frac{R}{2}}|\nabla u|^2\dx}\lesssim
\sup_{F}{\en{ |F u|^2}},
\end{equation}
where the supremum runs over all linear functionals $F$ bounded in the sense of
\begin{equation}\label{b9}
|Fu|^2\le\int_{|x|\le R}|\nabla u|^2\dx,
\end{equation}
and where here and in the proof $\lesssim$ means $\le C$ with a generic $C=C(d,\lambda)$.
\end{lemma}

We note that the statement of Lemma \ref{L} is trivial for an ensemble $\langle\cdot\rangle$
supported on a single function. Hence Lemma \ref{L} expresses
some compactness of ensembles of $a$-harmonic functions.


\section{Proofs}

{\sc Proof of Theorem \ref{T}}.
Throughout the proof $\lesssim$ denotes $\le C$, where $C$ is a generic constant that
only depends on the dimension $d$, the ellipticity ratio $\lambda>0$, and the exponent $\alpha\in(0,1)$.
By a rescaling of space we may assume w.\ l.\ o.\ g.\ that $r_*=1$;
by homogeneity, we may w.\ l.\ o.\ g.\ assume $(\int|g|^2)^\frac{1}{2}= 1$.
We set for abbreviation $R:=\frac{1}{4}|x_0|\ge 1$.
We first list and motivate the main steps in the proof.

\medskip

{\bf Step 1}. In the first step, we upgrade the excess decay (\ref{1.10}) in the sense
that we replace the optimal slope $\xi$ on the l.\ h.\ s.\ by a fixed slope that does not
depend on $r$, but is the optimal slope for some radius $r_0$ of order one. 
More precisely, suppose that (\ref{T.1}) holds for
a point $x$, say $x=0$. Then there exists a radius $1\le r_0\lesssim 1$ such that
for all radii $R\ge r_0$ and $a$-harmonic functions $u$ in $\{|x|\le R\}$, the optimal
slope on scale $r_0$
\begin{equation}\nonumber
\xi:={\rm argmin}\int_{B_{r_0}}|\nabla u-\xi_i(e_i+\nabla\phi_i)|^2
\end{equation}
is such that for $r\in[r_0,R]$:
\begin{equation}\label{1.6}
\left(\int_{B_{r}}|\nabla u-\xi_i(e_i+\nabla\phi_i)|^2\right)^\frac{1}{2}
\lesssim(\frac{r}{R})^{\frac{d}{2}+\alpha}\left(\int_{B_R}|\nabla u|^2\right)^\frac{1}{2}.
\end{equation}
In addition, we have
\begin{equation}\label{5.7}
|\xi|\lesssim\left(\fint_{B_{r}}|\nabla u|^2\right)^\frac{1}{2}
\lesssim\left(\fint_{B_R}|\nabla u|^2\right)^\frac{1}{2}.
\end{equation}

Since $r_0 \lesssim 1$, the second inequality holds, possibly with a worse constant, for all $1 \le r \le R$. 

\medskip

{\bf Step 2}. We have the following estimates on $\nabla u$ and $\nabla v$
\begin{equation}\label{2.1}
\left(\int|\nabla u|^2\right)^\frac{1}{2}\lesssim 1,\quad
{\rm sup}_{|x|\ge 2}(|x|^d|\nabla v(x)|+|x|^{d+1}|\nabla^2v(x)|)\lesssim 1.
\end{equation}
Moreover, $u$ and $v$ have vanishing ``constant invariant''
\begin{equation}\label{2.4}
\int\nabla\eta_r\cdot a\nabla u=\int\nabla\eta_r\cdot a_h\nabla v=0
\end{equation}
and, thanks to the special form of the r.\ h.\ s.\ of (\ref{2.3}),
identical ``linear invariants'' (for $k=1,\ldots,d$)
\begin{equation}\label{2.5}
\int\nabla\eta_r\cdot((x_k+\phi_k)a\nabla u-ua(e_k+\nabla\phi_k))
=\int\nabla\eta_r\cdot(x_ka_h\nabla v-va_he_k),
\end{equation}
for all $r\ge 1$, and where $\eta_r(x)=\eta(\frac{x}{r})$ and $\eta$ 
is a cut-off function for $\{|x|\le 1\}$ in $\{|x|\le 2\}$. We speak of
invariants, since for two $a$-harmonic functions $u$ and $\tilde u$ (in our case
$\tilde u=1$ for the constant invariant and $\tilde u=x_k+\phi_k$ for the linear
invariant) defined in
$\{|x|\ge 1\}$, the value of the boundary integral 
$\int_{\partial\Omega}\nu\cdot(\tilde ua\nabla u-ua\nabla\tilde u)$ does not 
depend on the open set $\Omega$ provided the latter contains $B_1$.

\medskip

{\bf Step 3}. We consider the homogenization error
\begin{equation}\label{3.3}
w:=u-(v+\tilde\phi_i\partial_i v),
\end{equation}
where $\tilde\phi$ denotes the following blended version of the corrector:
\begin{equation}\label{3.2}
\tilde\phi:=(1-\eta)(\phi-\fint_{B_1(0)}\phi)+\eta(\phi-\fint_{B_1(x_0)}\phi),
\end{equation}
where $\eta$ is a cut-off function for $\{|x-x_0|\le R\}$ in $\{|x-x_0|\le 2R\}$.
Then we have
\begin{equation}\label{3.1}
-\nabla\cdot a\nabla w=\nabla\cdot h\quad\mbox{for}\;|x|\ge 2,
\end{equation}
where $\nabla w$ and $h$ satisfy the following estimates
\begin{equation}\label{3.9}
\left(\int_{|x|\ge 2}|\nabla w|^2\right)^\frac{1}{2}\lesssim 1
\end{equation}
and
%
\begin{align}\label{3.14}
\left(\int_{|x|\ge r}|h|^2\right)^\frac{1}{2}&\lesssim\frac{1}{r^{\frac{d}{2}+\alpha}}\qquad
\textrm{for}\;r\ge 2,\\
\label{3.18}
\left(\int_{|x-x_0|\le r}|h|^2\right)^\frac{1}{2}&\lesssim\frac{r^\frac{d}{2}}{R^{d+\alpha}}\qquad
\textrm{for}\;1 \le r\le R.
\end{align}
Moreover, $w$ has an asymptotically vanishing constant invariant
\begin{equation}\label{3.19}
\lim_{r\uparrow\infty}\int\nabla\eta_r\cdot a\nabla w=0
\end{equation}
and asymptotically vanishing linear invariants
\begin{equation}\label{3.20}
\lim_{r\uparrow\infty}\int\nabla\eta_r\cdot ((x_k+\phi_k)a\nabla w-wa(e_k+\nabla\phi_k))=0
\quad\mbox{for}\;k=1,\ldots,d.
\end{equation}

\medskip

{\bf Step 4}. Extension into the origin. There exists $\bar w$ (with square integrable gradient), 
$\bar h$, and $\bar f$ defined on all of $\mathbb{R}^d$
such that 
\begin{equation}\label{4.5}
-\nabla\cdot a\nabla \bar w=\nabla\cdot\bar h+\bar f,
\end{equation}
while
\begin{equation}\label{4.6}
\nabla\bar w=\nabla w\quad\mbox{in}\;\{|x|\ge 4\},
\quad \supp\bar h\subset\{|x|\ge 2\},\quad \supp\bar f\subset\{|x|\le 4\}.
\end{equation}
The r.\ h.\ s.\ $\bar f$ and $\bar h$ satisfy the estimates
\begin{align}\label{4.1}
\left(\int\bar f^2\right)^\frac{1}{2}&\lesssim 1,\\
%
\intertext{and}
%
\label{4.2}
\left(\int_{|x|\ge r}|\bar h|^2\right)^\frac{1}{2}&\lesssim\frac{1}{r^{\frac{d}{2}+\alpha}}\quad
\mbox{for}\;r\ge 2,\\
%
%
\label{4.3}
\left(\int_{|x-x_0|\le r}|\bar h|^2\right)^\frac{1}{2}&\lesssim\frac{r^\frac{d}{2}}{R^{d+\alpha}}\quad
\mbox{for}\;1 \le r\le R.
\end{align}
Moreover, we have vanishing constant and linear invariants:
\begin{equation}\label{4.4}
\int \bar f=0,\quad\int((e_k+\nabla\phi_k)\cdot\bar h-(x_k+\phi_k)\bar f)=0\quad
\mbox{for}\;k=1,\ldots,d.
\end{equation}
The estimates (\ref{4.1})-(\ref{4.3}) ensure that the integrals in (\ref{4.4}) 
converge absolutely.

\medskip

{\bf Step 5}. From the equation (\ref{4.5}), the vanishing invariants (\ref{4.4}), 
and the estimates (\ref{4.1}) \& (\ref{4.2}) it follows
\begin{equation}\label{5.9}
\left(\int_{|x|\ge R}|\nabla\bar w| ^2\right)^\frac{1}{2}\lesssim\frac{\ln R}{R^{\frac{d}{2}+\alpha}}.
\end{equation}

\medskip

{\bf Step 6}. From the same ingredients as in Step 5 and in addition (\ref{4.3}),
one obtains the following localized version of (\ref{5.9}):
\begin{equation}\nonumber
\left(\int_{|x-x_0|\le 1}|\nabla\bar w|^2\right)^\frac{1}{2}\lesssim\frac{\ln R}{R^{d+\alpha}}.
\end{equation}

\medskip

{\bf Step 7}. Conclusion.


{\sc Argument for Step 1}

We first argue that there exists a radius $1\le r_0\sim 1$ such that
\begin{equation}\label{1.1}
\fint_{B_r}|\xi_i\cdot(e_i+\nabla\phi_i)|^2\sim|\xi|^2\quad\mbox{for}\;r\ge r_0.
\end{equation}
Indeed, (\ref{1.1}) easily follows from the sublinear growth (\ref{T.1})
of $\phi$ in the point $0$.
The upper bound is a consequence of Caccioppoli's estimate for the
$a$-harmonic function $x_i+(\phi_i-\fint_{B_{2r}}\phi_i)$, cf.\ (\ref{i2}), leading to
\begin{align*}
\left(\int_{B_r}|\xi_i\cdot(e_i+\nabla\phi_i)|^2\right)^\frac{1}{2}
&\lesssim\frac{1}{r}\left(\int_{B_{2r}}|\xi_i(x_i+(\phi_i-\fint_{B_{2r}}\phi_i))|^2\right)^\frac{1}{2}\\
&\lesssim|\xi|\left(1+\frac{1}{r}\left(\int_{B_{2r}}|\phi-\fint_{B_{2r}}\phi|^2\right)^\frac{1}{2}\right).
\end{align*}
The lower bound is a consequence of Jensen's inequality once we introduce
a cut-off function $\eta$ of $B_\frac{1}{2}$ in $B_1$ and set $\eta_r(x)=\eta(\frac{x}{r})$:
\begin{align*}
\left(\int_{B_r}|\xi_i\cdot(e_i+\nabla\phi_i)|^2\right)^\frac{1}{2}
&\ge\left(\int\eta_r|\xi_i\cdot(e_i+\nabla\phi_i)|^2\right)^\frac{1}{2}\\
&\ge{\textstyle\frac{1}{(\int\eta_r)^\frac{1}{2}}}\left|\int\eta_r\xi_i\cdot(e_i+\nabla\phi_i)\right|\\
&={\textstyle\frac{1}{(\int\eta_r)^\frac{1}{2}}}\left|\xi\int\eta_r-\xi_i\int(\phi_i-\fint_{B_r}\phi_i)\nabla\eta_r\right|\\
&\ge|\xi|\left(\frac{1}{C}r^\frac{d}{2}-Cr^{\frac{d}{2}-1}\left(\int_{B_r}|\phi_i-\fint_{B_r}\phi_i|^2\right)^\frac{1}{2}\right).
\end{align*}

\medskip

For an $a$-harmonic function $u$ in $B_R$ and a radius $r_0\le r\le R$ we consider
\begin{equation}\label{1.5}
\xi_r:={\rm argmin}\fint_{B_r}|\nabla u-\xi_i(e_i+\nabla\phi_i)|^2
\end{equation}
and claim that
\begin{equation}\label{1.7}
|\xi_r-\xi_R|\lesssim\left(\inf_{\xi}\fint_{B_R}|\nabla u-\xi_{i}(e_i+\nabla\phi_i)|^2\right)^\frac{1}{2}.
\end{equation}
Indeed, by the triangle inequality in $\mathbb{R}^d$ and since $\alpha>0$, it is enough to show
for $r_0\le r\le r'\le R$ with $r'\le 2r$ that
\begin{equation}\nonumber
|\xi_r-\xi_{r'}|\lesssim\p{\frac{r}{R}}^\alpha
\left(\inf_{\xi}\fint_{B_R}|\nabla u-\xi_{i}(e_i+\nabla\phi_i)|^2\right)^\frac{1}{2},
\end{equation}
which thanks to Theorem \ref{T2} follows from
\begin{equation}\nonumber
|\xi_r-\xi_{r'}|^2\le\inf_{\xi}\fint_{B_{r}}|\nabla u-\xi_{i}(e_i+\nabla\phi_i)|^2+
\inf_{\xi}\fint_{B_{r'}}|\nabla u-\xi_{i}(e_i+\nabla\phi_i)|^2,
\end{equation}
which by definition (\ref{1.5}) in turn follows from
\begin{equation}\nonumber
|\xi_r-\xi_{r'}|^2\le\fint_{B_{r}}|\nabla u-\xi_{r,i}(e_i+\nabla\phi_i)|^2+
\fint_{B_{r}}|\nabla u-\xi_{r',i}(e_i+\nabla\phi_i)|^2.
\end{equation}
The latter finally follows from the triangle inequality in $L^2(B_r)$ and from
(\ref{1.1}) applied to $\xi=\xi_r-\xi_{r'}$.

\medskip

Now (\ref{1.6}) is an easy consequence of (\ref{1.7}) (with $(r,R)$ replaced by $(r_0,r)$)
and (\ref{1.10}):
%
\begin{equation*}
\left(\int_{B_{r}}|\nabla u-\xi_i(e_i+\nabla\phi_i)|^2\right)^\frac{1}{2}
\stackrel{(\ref{1.7})}{\lesssim}
\left(\inf_{\tilde\xi}\int_{B_{r}}|\nabla u-\tilde\xi_i(e_i+\nabla\phi_i)|^2\right)^\frac{1}{2}
\stackrel{(\ref{1.10})}{\lesssim}
(\frac{r}{R})^{\frac{d}{2}+\alpha}\left(\int_{B_R}|\nabla u|^2\right)^\frac{1}{2}.
\renewcommand\arraystretch{1} 
\end{equation*}

\medskip

We finally turn to the argument for (\ref{5.7}). For this purpose we first note that
for all radii $r\ge r_0$
\begin{equation}\label{1.20}
|\xi_r|+\left(\inf_{\xi}\fint_{B_r}|\nabla u-\xi_i(e_i+\nabla\phi_i)|^2\right)^\frac{1}{2}
\lesssim \left(\fint_{B_r}|\nabla u|^2\right)^\frac{1}{2}.
\end{equation}
Indeed, that the second l.\ h.\ s.\ term is dominated by the r.\ h.\ s.\ is obvious.
For the first l.\ h.\ s.\ term we note by (\ref{1.1}), the triangle inequality in
$L^2(B_r)$, and the definition (\ref{1.5}) of $\xi_r$ that
\begin{align*}
|\xi_r|\lesssim\left(\fint_{B_r}|\xi_{r,i}(e_i+\nabla\phi_i)|^2\right)^\frac{1}{2}
\le\left(\inf_{\xi}\fint_{B_r}|\nabla u-\xi_i(e_i+\nabla\phi_i)|^2\right)^\frac{1}{2}
+\left(\fint_{B_r}|\nabla u|^2\right)^\frac{1}{2}.
\end{align*}

\medskip

Equipped with (\ref{1.20}), and more importantly (\ref{1.7}) and (\ref{1.10}), 
we may now tackle (\ref{5.7}). For the first
estimate in (\ref{5.7}), we appeal to the triangle inequality, (\ref{1.7}) with $(r,R)$ replaced
by $(r_0,r)$, and (\ref{1.20}):
\begin{align*}
|\xi_{r_0}|\le|\xi_{r}|+|\xi_{r_0}-\xi_r|
\stackrel{(\ref{1.7})}{\lesssim}
|\xi_{r}|+\left(\inf_{\xi}\fint_{B_r}|\nabla u-\xi_i(e_i+\nabla\phi_i)|^2\right)^\frac{1}{2}
\stackrel{(\ref{1.20})}{\lesssim}\left(\fint_{B_r}|\nabla u|^2\right)^\frac{1}{2}.
\end{align*}
For the second estimate in (\ref{5.7}), we use the triangle inequality in $L^2(B_r)$ and
in $\mathbb{R}^d$, (\ref{1.7}), (\ref{1.10}), and (\ref{1.20}):
%
\begin{equation*}
\renewcommand\arraystretch{1.8} 
\begin{array}{@{} r@{}c@{}l @{}}
 \left(\fint_{B_r}|\nabla u|^2\right)^\frac{1}{2}
&\stackrel{(\ref{1.5})}{\le}&
|\xi_{R}|+|\xi_{R}-\xi_r|+\left(\inf_{\xi}\fint_{B_r}|\nabla u-\xi_i(e_i+\nabla\phi_i)|^2\right)^\frac{1}{2}\nonumber\\
&\stackrel{(\ref{1.7})}{\lesssim}&
|\xi_{R}|+\max_{\rho=r,R}\left(\inf_{\xi}\fint_{B_\rho}|\nabla u-\xi_i(e_i+\nabla\phi_i)|^2\right)^\frac{1}{2}\nonumber\\
&\stackrel{(\ref{1.10})}{\lesssim}&
|\xi_{R}|+\left(\inf_{\xi}\fint_{B_R}|\nabla u-\xi_i(e_i+\nabla\phi_i)|^2\right)^\frac{1}{2}\nonumber\\
&\stackrel{(\ref{1.20})}{\lesssim}&\left(\fint_{B_R}|\nabla u|^2\right)^\frac{1}{2}.
\end{array}
\end{equation*}
%

\medskip

{\sc Argument for Step 2}. The first estimate in (\ref{2.1}) is an immediate consequence
of the definition (\ref{2.2}) and the energy estimate. We now turn to the second estimate
in (\ref{2.1}). We first note that the r.\ h.\ s.\ $\tilde g:=g_i(e_i+\partial_i\phi)$
of (\ref{2.3}) can be bounded with help of Caccioppoli's estimate
\begin{equation}\nonumber
\int|\tilde g|\le\left(\int g_i^2\right)^\frac{1}{2}\left(\int_{B_1}|e_i+\nabla\phi_i|^2\right)^\frac{1}{2}
\lesssim\left(\int_{B_2}|x+(\phi-\fint_{B_2}\phi)|^2\right)^\frac{1}{2} \lesssim 1.
\end{equation}
Since moreover, it is supported in $B_1$, we obtain the following representation for $|x|\ge 2$
\begin{equation}\nonumber
\nabla v(x)=\int\nabla^2 G_h(x-y)\tilde g(y)\dy,\quad
\nabla^2 v(x)=\int\nabla^3 G_h(x-y)\tilde g(y)\dy
\end{equation}
in terms of the constant-coefficient Green's function $G_h$. The second estimate in
(\ref{2.1}) now follows using the homogeneity of $G_h$.

\medskip

We now turn to the invariants. The first identity (\ref{2.4}) follows immediately from 
integration by parts of the equations (\ref{2.2}) \& (\ref{2.3}) using the fact that
the respective right hand sides are supported in $B_1$. We now turn to the second identity (\ref{2.5}),
which also follows from integration by parts, but this time relying on the Green formulas
\begin{align*}
\nabla\cdot[(x_k+\phi_k)(a\nabla u+g)-ua(e_k+ \nabla\phi_k)]&=(e_k+\nabla\phi_k)\cdot g,\\
\nabla\cdot[x_k(a_h\nabla v+g_i(e_i+\partial_i\phi))- v a_he_k]&=e_k\cdot [g_i(e_i+\partial_i\phi)],
\end{align*}
which follow from (\ref{2.2}) \& (\ref{2.3}) in conjunction with (\ref{i2}),
and the pointwise identity
\begin{equation}\nonumber
(e_k+\nabla\phi_k)\cdot g=e_k\cdot[g_i(e_i+\partial_i\phi)].
\end{equation}

\medskip

{\sc Argument for Step 3}. We start by establishing the formula (\ref{3.1}) with
\begin{align}\label{3.15}
h&:=(\tilde \phi_i a-\tilde\sigma_i)\nabla\partial_iv
+\partial_iv\Big((\fint_{B_1(0)}\phi_i-\fint_{B_1(x_0)}\phi_i) a-(\fint_{B_1(0)}\sigma_i-\fint_{B_1(x_0)}\sigma_i)\Big)\nabla\eta,
\end{align}
where, in line with (\ref{3.2}), we have set
\begin{equation}\label{3.5}
\tilde\sigma:=(1-\eta)(\sigma-\fint_{B_1(0)}\sigma)+\eta(\sigma-\fint_{B_1(x_0)}\sigma).
\end{equation}
Indeed, from definition (\ref{3.3}) \& (\ref{3.2}) we obtain
\begin{equation}\label{3.10}
\nabla w=\nabla u-[\partial_iv(e_i+\nabla\phi_i+(\fint_{B_1(0)}\phi_i-\fint_{B_1(x_0)}\phi_i)\nabla\eta)
+\tilde\phi_i\nabla\partial_iv]
\end{equation}
and thus, using $\nabla\cdot a(e_i+\nabla\phi_i)=0$, cf.\ (\ref{i2}), 
\begin{align*}
-\nabla\cdot a\nabla w=-\nabla\cdot a\nabla u+\nabla\partial_iv \cdot a(e_i+\nabla\phi_i)
+\nabla\cdot[\partial_iv(\fint_{B_1(0)}\phi_i-\fint_{B_1(x_0)}\phi_i) a\nabla\eta
+\tilde\phi_ia\nabla\partial_iv].
\end{align*}
Using the identity $a(e_i+\nabla\phi_i)=a_he_i+\nabla\cdot\sigma_i$, cf.\ (\ref{i5}), we have
$\nabla\partial_iv \cdot a(e_i+\nabla\phi_i)=\nabla\cdot a_h\nabla v+\nabla\partial_iv\cdot(\nabla\cdot\sigma_i)$.
Using that the r.\ h.\ s.\ of (\ref{2.2}) \& (\ref{2.3}) are supported in $B_1$, we thus obtain in $\{|x|\ge 1\}$
\begin{equation}\nonumber
-\nabla\cdot a\nabla w=\nabla\partial_iv\cdot(\nabla\cdot\sigma_i)
+\nabla\cdot[\partial_iv(\fint_{B_1(0)}\phi_i-\fint_{B_1(x_0)}\phi_i) a\nabla\eta
+\tilde\phi_ia\nabla\partial_iv].
\end{equation}
It remains to substitute $\sigma$ by $\tilde\sigma$ and to bring the related part of the above r.\ h.\ s.\ into
divergence form. Indeed, by definition (\ref{3.5})
\begin{equation}\nonumber
\nabla\cdot\sigma_i
=\nabla\cdot\tilde\sigma_i-\nabla\cdot(\eta(\fint_{B_1(0)}\sigma_i-\fint_{B_1(x_0)}\sigma_i))
=\nabla\cdot\tilde\sigma_i-(\fint_{B_1(0)}\sigma_i-\fint_{B_1(x_0)}\sigma_i)\nabla\eta,
\end{equation}
and by the identities
\begin{equation}\label{3.21}
\nabla\zeta\cdot(\nabla\cdot\sigma)=\nabla\cdot(\zeta\nabla\cdot\sigma)
=-\nabla\cdot(\sigma\nabla\zeta)\quad\mbox{for skew}\;\sigma,
\end{equation}
we obtain as desired
\begin{equation}\nonumber
\nabla\partial_iv\cdot(\nabla\cdot\sigma_i)
=-\nabla\cdot[\tilde\sigma_i\nabla\partial_iv]
-\nabla\cdot[\partial_iv(\fint_{B_1(0)}\sigma_i-\fint_{B_1(x_0)}\sigma_i)\nabla\eta].
\end{equation}

\medskip

In order to estimate the contribution to $\nabla w$ and $h$ that comes from the difference
$\fint_{B_1(0)}-\fint_{B_1(x_0)}$ of the average values we now argue that
\begin{equation}\label{3.6}
\left|\fint_{B_1(0)}(\phi,\sigma)-\fint_{B_1(x_0)}(\phi,\sigma)\right|\lesssim R^{1-\alpha}.
\end{equation}
To keep notation light, we write $\phi$ instead of $(\phi,\sigma)$.
Let us first argue how to reduce (\ref{3.6}) to
\begin{equation}\label{3.7}
\left|\fint_{B_r(x)}\phi-\fint_{B_1(x)}\phi\right|\lesssim r^{1-\alpha}
\quad\mbox{for}\;r\ge 1,\;x\in\{0,x_0\}.
\end{equation}
Indeed, (\ref{3.6}) follows from (\ref{3.7}) and our assumption (\ref{T.1}) via
the string of inequalities
\begin{eqnarray*}
\lefteqn{\left|\fint_{B_1(0)}\phi-\fint_{B_1(x_0)}\phi\right|}\\
&\le&\left(\fint_{B_{2R}(\frac{1}{2}x_0)}|\phi-\fint_{B_1(0)}\phi|^2\right)^\frac{1}{2}
+\left(\fint_{B_{2R}(\frac{1}{2}x_0)}|\phi-\fint_{B_1(x_0)}\phi|^2\right)^\frac{1}{2}\\
&\lesssim&\left(\fint_{B_{4R}(0)}|\phi-\fint_{B_1(0)}\phi|^2\right)^\frac{1}{2}
+\left(\fint_{B_{4R}(x_0)}|\phi-\fint_{B_1(x_0)}\phi|^2\right)^\frac{1}{2}\\
&\lesssim&\left(\fint_{B_{4R}(0)}|\phi-\fint_{B_{4R}(0)}\phi|^2\right)^\frac{1}{2}
+\left(\fint_{B_{4R}(x_0)}|\phi-\fint_{B_{4R}(x_0)}\phi|^2\right)^\frac{1}{2}\\
&&+\left|\fint_{B_{4R}(0)}\phi-\fint_{B_1(0)}\phi\right|+
\left|\fint_{B_{4R}(x_0)}\phi-\fint_{B_1(x_0)}\phi\right| \overset{\eqref{T.1},\eqref{3.7}}{\lesssim} R^{1-\alpha}.
\end{eqnarray*}
For (\ref{3.7}), we focus on $x=0$ and note that by a decomposition into
dyadic radii, it is enough to show 
\begin{equation}\nonumber
\left|\fint_{B_r}\phi-\fint_{B_{r'}}\phi\right|\lesssim r^{1-\alpha}
\quad\mbox{for}\;2r\ge r'\ge r\ge 1.
\end{equation}
This estimate follows from a similar string of inequalities as the one before:
\begin{align*}
\left|\fint_{B_r}\phi-\fint_{B_{r'}}\phi\right|
&\le\left(\fint_{B_{r}}|\phi-\fint_{B_r}\phi|^2\right)^\frac{1}{2}
+\left(\fint_{B_{r}}|\phi-\fint_{B_{r'}}\phi|^2\right)^\frac{1}{2}\\
&\lesssim\left(\fint_{B_r}|\phi-\fint_{B_r}\phi|^2\right)^\frac{1}{2}
+\left(\fint_{B_{r'}}|\phi-\fint_{B_{r'}}\phi|^2\right)^\frac{1}{2}
\end{align*}
and an application of (\ref{T.1}).

\medskip

We now turn to (\ref{3.9}). We start from the formula (\ref{3.10}), which we rewrite as
\begin{align*}
\nabla w&=\nabla u-\partial_iv\Big(e_i+\nabla\phi_i+(\fint_{B_1(0)}\phi_i-\fint_{B_1(x_0)}\phi_i)\nabla\eta\Big)\\
&\quad-\Big((\phi_i-\fint_{B_1(0)}\phi_i)+\eta(\fint_{B_1(0)}\phi_i-\fint_{B_1(x_0)}\phi_i)\Big)\nabla\partial_iv.
\end{align*}
From the estimate (\ref{3.6}) on averages and the estimate (\ref{2.1}) on $v$
we thus obtain for $|x|\ge 2$
\begin{align}\label{3.13}
|\nabla w|\lesssim|\nabla u|+|x|^{-d}\big(|{\rm id}+\nabla\phi|+R^{1-\alpha}|\nabla\eta|\big)
+|x|^{-(d+1)}\big(|\phi-\fint_{B_{1}(0)}\phi|+R^{1-\alpha}\eta\big).
\end{align}
Taking the $L^2(\{|x|\ge 2\})$-norm, we see that the $\nabla u$-term is bounded according to (\ref{2.1}).
By the choice of the cut-off function $\eta$, the function $|x|^{-d}|\nabla\eta|+|x|^{-(d+1)}\eta$
is supported in $\{|x-x_0|\le 2R\}$ and bounded by $R^{-(d+1)}$ so that the contribution of
this term to $\left(\int_{|x|\ge 2}|\nabla w|^2\right)^\frac{1}{2}$ is estimated by
$R^{-\alpha-\frac{d}{2}}$. We turn to the term involving $\phi$ and will
for later purposes show the slightly more general statement for $r\ge 1$
\begin{equation}\label{3.17}
\left(\int_{|x|\ge r}(|x|^{-(d+1)}|(\phi,\sigma)-\fint_{B_{1}(0)}(\phi,\sigma)|)^2\right)^\frac{1}{2}
\lesssim\frac{1}{r^{\frac{d}{2}+\alpha}},
\end{equation}
which, restricting to $\phi$ in our notation, follows by dyadic summation from 
\begin{equation}\nonumber
\left(\int_{r\le |x|\le 2r}|\phi-\fint_{B_{1}(0)}\phi|^2\right)^\frac{1}{2}
\lesssim r^{\frac{d}{2}+1-\alpha},
\end{equation}
which trivially follows from 
\begin{equation}\nonumber
\left(\fint_{B_r(0)}|\phi-\fint_{B_{1}(0)}\phi|^2\right)^\frac{1}{2}
\lesssim r^{1-\alpha}.
\end{equation}
The last estimate is a combination of our assumption (\ref{T.1}) with (\ref{3.7}).
We now turn to the estimate of the $\nabla\phi$-term in (\ref{3.13}). By dyadic summation, it is
enough to show
\begin{equation}\nonumber
\left(\int_{r\le |x|\le 2r}|{\rm id}+\nabla\phi|^2\right)^\frac{1}{2}
\lesssim r^{\frac{d}{2}},
\end{equation}
which follows from Caccioppoli's estimate and assumption (\ref{T.1}):
\begin{align*}
\left(\fint_{B_{2r}}|{\rm id}+\nabla\phi|^2\right)^\frac{1}{2}
\lesssim\frac{1}{r}\left(\fint_{B_{4r}}|x+(\phi-\fint_{B_{4r}}\phi)|^2\right)^\frac{1}{2}
\lesssim 1+\frac{1}{r}\left(\fint_{B_{4r}}|\phi-\fint_{B_{4r}}\phi|^2\right)^\frac{1}{2}.
\end{align*}

\medskip

We now turn to estimate (\ref{3.14}). For this purpose, we rewrite the definition (\ref{3.15}) of $h$ as
\begin{align}
h&=\Big((\phi_i-\fint_{B_1(0)}\phi_i) a-(\sigma_i-\fint_{B_1(0)}\sigma_i)\Big)\nabla\partial_iv\nonumber\\
&\quad+\eta\Big((\fint_{B_1(0)}\phi_i-\fint_{B_1(x_0)}\phi_i) a-(\fint_{B_1(0)}\sigma_i-\fint_{B_1(x_0)}\sigma_i)\Big)
\nabla\partial_iv\nonumber\\
&\quad+\partial_iv\Big((\fint_{B_1(0)}\phi_i-\fint_{B_1(x_0)}\phi_i) a-(\fint_{B_1(0)}\sigma_i-\fint_{B_1(x_0)}\sigma_i)\Big)\nabla\eta.
\nonumber
\end{align}
Inserting, as above, the estimate (\ref{2.1}) on $v$ and the estimate (\ref{3.6}) on the averages we obtain
for $|x|\ge 2$
\begin{equation}\label{3.16}
|h|\lesssim|x|^{-(d+1)}\Big|(\phi,\sigma)-\fint_{B_1(0)}(\phi,\sigma)\Big|
+R^{1-\alpha}\big(|x|^{-(d+1)}\eta+|x|^{-d}|\nabla\eta|\big).
\end{equation}
Since by definition of $\eta$, the function $|x|^{-(d+1)}\eta+|x|^{-d}|\nabla\eta|$ 
is supported in $\{|x|\ge 2R\}$ and bounded by $R^{-(d+1)}$, 
its contribution to $\left(\int_{|x|\ge r}|h|^2\right)^\frac{1}{2}$
vanishes for $r\le R$ and is bounded by $R^\frac{d}{2}R^{-(d+\alpha)}\le r^{-(\frac{d}{2}+\alpha)}$
for $r\ge R$. The first r.\ h.\ s.\ term in (\ref{3.16}) was treated in (\ref{3.17}).

\medskip

We finally turn to the last estimate (\ref{3.18}). For this purpose, we note that
because of the choice of $\eta$, on $B_R(x_0)$, the definition (\ref{3.15}) of $h$ turns into
\begin{equation*}
h=\Big((\phi_i-\fint_{B_1(x_0)}\phi_i) a-(\sigma_i-\fint_{B_1(x_0)}\sigma_i)\Big)\nabla\partial_iv,
\end{equation*}
so that by the estimate (\ref{2.1}) on $v$ we have for $1 \le r\le R$
\begin{equation*}
\left(\fint_{B_r(x_0)}|h|^2\right)^\frac{1}{2}\lesssim R^{-(d+1)}
\left(\fint_{B_r(x_0)}|(\phi,\sigma)-\fint_{B_1(x_0)}(\phi,\sigma)|^2\right)^\frac{1}{2},
\end{equation*}
which by estimate (\ref{3.7}) on averages turns into
\begin{equation*}
\left(\fint_{B_r(x_0)}|h|^2\right)^\frac{1}{2}\lesssim R^{-(d+1)}\left(
\left(\fint_{B_r(x_0)}|(\phi,\sigma)-\fint_{B_r(x_0)}(\phi,\sigma)|^2\right)^\frac{1}{2}+r^{1-\alpha}\right),
\end{equation*}
so that the desired estimate in the strengthened form of
\begin{equation*}
\left(\int_{B_r(x_0)}|h|^2\right)^\frac{1}{2} \lesssim \frac{r^{\frac{d}{2}+1-\alpha}}{R^{d+1}}
\end{equation*}
now follows from assumption (\ref{T.1}).

\medskip

We finally turn to the asymptotic invariants (\ref{3.19}) \& (\ref{3.20}). Since we may assume
$r\ge 6R$, we may ignore the presence of $\eta$ in the definition (\ref{3.3}) of $w$ 
(and assume w.\ l.\ o.\ g.\ that $\int_{B_1(0)}(\phi,\sigma)=0$ for notational simplicity).
Hence we have
\begin{equation}\label{3.23}
a\nabla w\stackrel{\phantom{(\ref{i5})}}{=}a\nabla u-\big(\partial_iva(e_i+\nabla\phi_i)+\phi_ia\nabla\partial_iv\big)
\stackrel{(\ref{i5})}{=}a\nabla u-\big(a_h\nabla v+\partial_iv\nabla\cdot\sigma_i+\phi_ia\nabla\partial_iv\big).
\end{equation}
Using the formula 
\begin{equation}\label{3.30}
\int\nabla\eta_r\cdot(\zeta\nabla\cdot\sigma+\sigma\nabla\zeta)=0\quad\mbox{for skew}\;\sigma,
\end{equation}
which follows from the last identity in (\ref{3.21}) (and even holds if $\zeta$ a priori
is not defined in $\{|x|\le r\}$ since it can be arbitrarily extended), we derive
\begin{equation}\label{3.24}
\int\nabla\eta_r\cdot a\nabla w
=\int\nabla\eta_r\cdot\big(a\nabla u-a_h\nabla v-(\phi_ia-\sigma_i)\nabla\partial_iv\big),
\end{equation}
which by the identity (\ref{2.4}) of the variants for $u$ and $v$ collapses into
\begin{equation}\nonumber
\int\nabla\eta_r\cdot a\nabla w
=-\int\nabla\eta_r\cdot\big((\phi_ia-\sigma_i)\nabla\partial_iv\big).
\end{equation}
Using the estimates (\ref{2.1}) on $v$, this yields
\begin{equation}\nonumber
\left|\int\nabla\eta_r\cdot a\nabla w\right|\lesssim r^{-(d+2)}\int_{B_{2r}}|(\phi,\sigma)|
\lesssim r^{-2}\left(\fint_{B_{2r}}|(\phi,\sigma)|^2\right)^\frac{1}{2}.
\end{equation}
Together with our assumption (\ref{T.1}), this implies (\ref{3.19}).

\medskip

We conclude this step with the argument for (\ref{3.20}). From the identity (\ref{3.23}) we deduce
\begin{align*}
(x_k+&\phi_k)a\nabla w-wa(e_k+\nabla\phi_k)\\
&=\big((x_k+\phi_k)a\nabla u-u a(e_k+\nabla\phi_k)\big)-\big( x_ka_hv-va_he_k\big)-\phi_ka_h\nabla v\\
&\quad-(x_k+\phi_k)\big(\partial_iv\nabla\cdot\sigma_i+\phi_ia\nabla\partial_iv\big)
+\phi_i\partial_iva(e_k+\nabla\phi_k)+v\nabla\cdot\sigma_k.
\end{align*}
Using that the linear invariants of $u$ and $v$ coincide, cf.\ (\ref{2.5}), and the 
formula (\ref{3.30}), only the following terms survive after application of $\int\nabla\eta_r\cdot$:
\begin{align*}
-\phi_ka_h\nabla v+\sigma_i\nabla\big((x_k+\phi_k)\partial_iv\big)
-(x_k+\phi_k)\phi_ia\nabla\partial_iv+\phi_i\partial_iva(e_k+\nabla\phi_k)
-\sigma_k \nabla v\\
=\partial_iv(\phi_ia+\sigma_i)(e_k+\nabla\phi_k)-(\phi_ka_h+\sigma_k)\nabla v
-(x_k+\phi_k)(\phi_ia-\sigma_i)\nabla\partial_iv.
\end{align*}
Together with the estimates (\ref{2.1}) on $v$, this implies by Caccioppoli's estimate
\begin{align*}
\bigg|\int&\nabla\eta_r\cdot\big((x_k+\phi_k)a\nabla w-wa(e_k+\nabla\phi_k)\big)\bigg|\\
&\lesssim
r^{-(d+1)}\int_{B_{2r}}|(\phi,\sigma)|(1+|{\rm id}+\nabla\phi|)
+r^{-(d+2)}\int_{B_{2r}}|(\phi,\sigma)||x+\phi|\\
&\lesssim
\frac{1}{r}\left(\fint_{B_{2r}}|(\phi,\sigma)|^2\right)^\frac{1}{2}
\left(1+\left(\fint_{B_{2r}}|{\rm id}+\nabla\phi|^2\right)^\frac{1}{2}
+\frac{1}{r}\left(\fint_{B_{2r}}|x+\phi|^2\right)^\frac{1}{2}\right).\\
&\lesssim
\frac{1}{r}\left(\fint_{B_{2r}}|(\phi,\sigma)|^2\right)^\frac{1}{2}
\left(1+
\frac{1}{r}\left(\fint_{B_{4r}}|\phi|^2\right)^\frac{1}{2}\right).
\end{align*}
By assumption (\ref{T.1}), this yields (\ref{3.20}).

\medskip

{\sc Argument for Step 4}. Select a cut-off function $\eta$ for $\{|x|\ge 4\}$ in $\{|x|\ge 2\}$
and set
\begin{equation}\nonumber
\bar w:=\eta(w-w_0),\quad \bar h:=\eta h-(w-w_0)a\nabla\eta,\quad
\bar f:=-\nabla\eta\cdot(a\nabla w+h),
\end{equation}
where $w_0$ is the average of $w$ on the annulus $\{2\le |x|\le 4\}$. By the choice
of $\eta$, (\ref{4.6}) is clearly satisfied. Since
\begin{equation}\nonumber
\nabla\bar w=\eta \nabla w+(w-w_0)\nabla\eta
\quad\mbox{and thus}\quad
-a\nabla \bar w-\bar h=-\eta (a\nabla w+h),
\end{equation}
we learn from (\ref{3.1}) that also (\ref{4.5}) holds. The estimate (\ref{4.1}) on $\bar f$ follows
from the estimates (\ref{3.9}) and (\ref{3.14}) on $w$ and $h$. As for estimate
(\ref{4.2}) on $\bar h$, we note that for the contribution $\eta h$ to $\bar h$,
the estimate immediately translates from (\ref{3.14}); for the contribution $-(w-w_0)a\nabla\eta$
we note that it is supported in $\{2\le |x|\le 4\}$ and estimated by $|w-w_0|$, and thus
the desired estimate follows from Poincar\'e's inequality with mean value zero on the
annulus $\{2\le|x|\le 4\}$. For estimate (\ref{4.3}) on $\bar h$, we note that on 
$B_R(x_0)$, $\bar h$ coincides with $h$, so that it follows immediately from (\ref{3.18}).

\medskip

We now turn to the invariants (\ref{4.4}). For the constant invariant we note that
we obtain from the equation (\ref{4.5}) for $r\ge 2$
\begin{equation}\nonumber
\int\bar f\stackrel{(\ref{4.6})}{=}\int\eta_r\bar f=\int\nabla\eta_r\cdot(a\nabla\bar w+\bar h)
\stackrel{(\ref{4.6})}{=}\int\nabla\eta_r\cdot(a\nabla w+ h).
\end{equation}
Hence the first identity in (\ref{4.4}) follows from (\ref{3.19}) provided we have
$\lim_{r\uparrow\infty}\int\nabla\eta_r\cdot h=0$. The latter is a consequence of (\ref{3.14}):
\begin{equation}\nonumber
\left|\int\nabla\eta_r\cdot h\right|\lesssim\frac{1}{r}\int_{r\le|x|\le2r}|h|
\le r^{\frac{d}{2}-1}\left(\int_{|x|\ge r}|h|^2\right)^\frac{1}{2}.
\end{equation}

\medskip

We now turn to the linear invariants in (\ref{4.4}). We start with the identity
which follows from the equations for $\bar w$ and $\phi_k$ and the support properties
\begin{align*}
\int\eta_r&\big(-(e_k+\nabla\phi_k)\cdot\bar h+(x_k+\phi_k)\bar f\big)\\
&\stackrel{(\ref{4.5})}{=}\int\nabla(\eta_r(x_k+\phi_k))\cdot a\nabla\bar w+(x_k+\phi_k)\nabla\eta_r\cdot\bar h\\
&\stackrel{\;(\ref{i2})\;}{=}\int\nabla\eta_r\cdot\big((x_k+\phi_k)(a\nabla\bar w+\bar h)-(\bar w-w_0)a(e_k+\nabla\phi_k)\big)\\
&\stackrel{(\ref{4.6})}{=}\int\nabla\eta_r\cdot\big((x_k+\phi_k)(a\nabla w+h)-wa(e_k+\nabla\phi_k)\big).
\end{align*}
Hence the linear invariants follow from (\ref{3.20}) once we show
\begin{equation}\label{4.8}
\int|-(e_k+\nabla\phi_k)\cdot\bar h+(x_k+\phi_k)\bar f|<\infty
\end{equation}
and
\begin{equation}\label{4.7}
\lim_{r\uparrow\infty}\int(x_k+\phi_k)\nabla\eta_r\cdot h=0.
\end{equation}
The limit (\ref{4.7}) follows from the estimate
\begin{align}
\left|\int(x_k+\phi_k)\nabla\eta_r\cdot h\right|
&\stackrel{\phantom{(\ref{3.14})}}{\lesssim}
\left(r^\frac{d}{2}+\frac{1}{r}\left(\int_{|x|\le 2r}|\phi_k|^2\right)^\frac{1}{2}\right)
\left(\int_{|x|\ge r}|h|^2\right)^\frac{1}{2}\nonumber\\
&\stackrel{(\ref{3.14})}{\lesssim}
\frac{1}{r^\alpha}\left(1+\frac{1}{r}\left(\fint_{B_{2r}}|\phi_k|^2\right)^\frac{1}{2}\right)\label{4.10}
\end{align}
together with the observation that the term $\frac{1}{r}\left(\fint_{B_{2r}}|\phi_k|^2\right)^\frac{1}{2}$
is of higher order since w.\ l.\ o.\ g.\ we may assume $\fint_{B_1}\phi_k=0$ so that
in view of (\ref{3.7}) we may appeal to our assumption (\ref{T.1}).

\medskip

We now turn to (\ref{4.8}); in view of the square integrability of $\bar f$ and $\bar h$
(established above) and the local square integrability of
$x_k+\phi_k$ and its gradient, it remains to show
\begin{equation}\label{4.9}
\int_{|x|\ge 2}|(e_k+\nabla\phi_k)\cdot h|<\infty.
\end{equation}
To this purpose, we divide into dyadic annuli and use Caccioppoli's estimate
\begin{align*}
\int_{r\le |x|\le 2r}|(e_k+\nabla\phi_k)\cdot h|&\le
\left(\int_{B_{2r}}|e_k+\nabla\phi_k|^2\right)^\frac{1}{2}
\left(\int_{|x|\ge r}|h|^2\right)^\frac{1}{2}\\
&\lesssim
\frac{1}{r}\left(\int_{B_{4r}}|x_k+\phi_k|^2\right)^\frac{1}{2}
\left(\int_{|x|\ge r}|h|^2\right)^\frac{1}{2}\\
&\lesssim
\frac{1}{r}\left(r^\frac{d}{2}+\left(\int_{B_{4r}}|\phi_k|^2\right)^\frac{1}{2}\right)
\left(\int_{|x|\ge r}|h|^2\right)^\frac{1}{2}.
\end{align*}
We now appeal to the same argument as for (\ref{4.10}) to see
\begin{align*}
\int_{r\le |x|\le 2r}|(e_k+\nabla\phi_k)\cdot h|\lesssim\frac{1}{r^\alpha}.
\end{align*}
Summation over dyadic annuli yields (\ref{4.9}).

\medskip

{\sc Argument for Step 5}. We give an argument by duality and therefore consider
for arbitrary square-integrable vector field $\tilde h$ supported in $\{|x|\ge R\}$
the finite energy solution of
\begin{equation}\label{5.6}
-\nabla\cdot a\nabla \tilde w=\nabla\cdot\tilde h.
\end{equation}
Since $\bar w$ is a finite energy solution of (\ref{4.5}) we have
\begin{equation}\label{5.1}
-\int\nabla\bar w\cdot\tilde h=-\int\nabla\bar w\cdot a\nabla\tilde w=\int(\nabla\tilde w\cdot \bar h
-\tilde w\bar f).
\end{equation}
Recall Step 1 and consider
\begin{equation}\label{5.5}
\xi:={\rm argmin}\int_{B_{r_0}}|\nabla\tilde w-\xi_k(e_k+\nabla\phi_k)|^2
\quad\mbox{and}\quad c:=\fint_{B_4}(\tilde w-\xi_k(x_k+\phi_k)).
\end{equation}
By the vanishing invariants (\ref{4.4}) we may post-process (\ref{5.1}) to
\begin{equation}\nonumber
-\int\nabla\bar w\cdot\tilde h=\int\Big(\big(\nabla\tilde w-\xi_k(e_k+\nabla\phi_k)\big)\cdot \bar h
-\big(\tilde w-\xi_k(x_k+\phi_k)-c\big)\bar f\Big).
\end{equation}
By the support conditions (\ref{4.6}) on $\bar f$ and $\bar h$, this implies
\begin{align*}
\bigg|&\int\nabla\bar w\cdot\tilde h\bigg|\\
&\le
\sum_{n=1}^\infty\int_{2^n\le|x|\le 2^{n+1}}|\nabla\tilde w-\xi_k(e_k+\nabla\phi_k)||\bar h|
+\int_{|x|\le 4}|\tilde w-\xi_k(x_k+\phi_k)-c||\bar f|\\
&\le
\sum_{n=1}^\infty\left(\int_{B_{2^{n+1}}}|\nabla\tilde w-\xi_k(e_k+\nabla\phi_k)\big|^2\right)^\frac{1}{2}
\left(\int_{|x|\ge 2^n}|\bar h|^2\right)^\frac{1}{2}\\
&\quad+\left(\int_{B_4}|\tilde w-\xi_k(x_k+\phi_k)-c|^2\right)^\frac{1}{2}\left(\int|\bar f|^2\right)^\frac{1}{2}.
\end{align*}
Inserting the estimates (\ref{4.1}) and (\ref{4.2}) on $\bar f$ and $\bar h$ and using
Poincar\'e's inequality on $B_4$ (in view of the definition (\ref{5.5}) of $c$)
\begin{equation}\label{5.8}
\left|\int\nabla\bar w\cdot\tilde h\right|\lesssim
\sum_{n=1}^\infty(\frac{1}{2^n})^{\alpha}\left(\fint_{B_{2^{n}}}
|\nabla\tilde w-\xi_k(e_k+\nabla\phi_k)\big|^2\right)^\frac{1}{2}.
\end{equation}
We now distinguish the cases of $2^n\le R$ and $2^n\ge R$. In case of $2^n\le R$,
since by assumption $\tilde h$ is supported in $\{|x|\ge R\}$, $\tilde w$ is $a$-harmonic in $B_R$, cf.\ (\ref{5.6}). We thus may appeal to (\ref{1.6}) in Step 1 and obtain in view of the definition of $\xi$ in (\ref{5.5})
\begin{equation}
\left(\fint_{B_{2^{n}}}|\nabla\tilde w-\xi_k(e_k+\nabla\phi_k)\big|^2\right)^\frac{1}{2}
\lesssim(\frac{2^n}{R})^\alpha\left(\fint_{B_{R}}|\nabla\tilde w|^2\right)^\frac{1}{2}.
\end{equation}
In case of $2^n\ge R$, we use (\ref{5.7}) in Step 1 (and once more Caccioppoli's estimate
in conjunction with assumption (\ref{T.1})) to obtain
\begin{align*}
\left(\fint_{B_{2^{n}}}|\nabla\tilde w-\xi_k(e_k+\nabla\phi_k)\big|^2\right)^\frac{1}{2}
&\lesssim\left(\fint_{B_{2^n}}|\nabla\tilde w|^2\right)^\frac{1}{2}
+|\xi|
\left(\fint_{B_{2^n}}|{\rm id}+\nabla\phi|^2\right)^\frac{1}{2}\\
&\lesssim\left(\fint_{B_{2^n}}|\nabla\tilde w|^2\right)^\frac{1}{2}
+\left(\fint_{B_{R}}|\nabla\tilde w|^2\right)^\frac{1}{2}.
\end{align*}
Using the energy inequality for (\ref{5.6}) in form of
\begin{equation}\nonumber
\left(\fint_{B_r}|\nabla\tilde w|^2\right)^\frac{1}{2}\lesssim
r^{-\frac{d}{2}}\left(\int|\nabla\tilde w|^2\right)^\frac{1}{2}
\lesssim r^{-\frac{d}{2}}\left(\int|\tilde h|^2\right)^\frac{1}{2},
\end{equation}
we obtain in either case
\begin{equation}\nonumber
\left(\fint_{B_{2^{n}}}|\nabla\tilde w-\xi_k(e_k+\nabla\phi_k)\big|^2\right)^\frac{1}{2}
\lesssim
       \left\{\begin{array}{ccc}(\frac{2^n}{R})^\alpha(\frac{1}{R})^{\frac{d}{2}}&\mbox{for}&2^n\le R\\
                                (\frac{1}{R})^{\frac{d}{2}}&\mbox{for}&2^n\ge R\end{array}\right\}
\left(\int|\tilde h|^2\right)^\frac{1}{2}.
\end{equation}
Inserting this into (\ref{5.8}) we obtain
\begin{equation*}
\left|\int\nabla\bar w\cdot\tilde h\right|\lesssim\frac{\ln R}{R^{\alpha+\frac{d}{2}}}\left(\int|\tilde h|^2\right)^\frac{1}{2}.
\end{equation*}
Since the only constraint on $\tilde h$ was that it is supported in $\{|x|\ge R\}$, we obtain (\ref{5.9}).

\medskip

{\sc Argument for Step 6}. In view of (\ref{5.9}) in Step 5, it is enough to show
\begin{equation}\label{6.6}
\left(\int_{B_1(x_0)}|\nabla\bar w|^2\right)^\frac{1}{2}\lesssim
\frac{\ln R}{R^{d+\alpha}}+\frac{1}{R^\frac{d}{2}}
\left(\int_{B_R(x_0)}|\nabla\bar w|^2\right)^\frac{1}{2},
\end{equation}
where w.\ l.\ o.\ g.\ we may assume that $R=2^N$ is dyadic.
To this purpose, we decompose the r.\ h.\ s.\ $\bar h$ of (\ref{4.5})
into $(\bar h_n)_{n=0,\ldots,N+1}$ such that $\bar h_0$ is supported in $B_1(x_0)$,
$\bar h_n$ is supported in the annulus $\{2^{n-1}\le|x-x_0|\le 2^n\}$ for $n=1,\ldots,N$,
and $\bar h_{N+1}$ is supported in the exterior domain $\{|x-x_0|\ge R\}$.
For $n=0,\ldots,N$ let $\bar w_n$ denote the finite energy solution of 
\begin{equation}\label{6.1}
-\nabla\cdot a\nabla \bar w_n=\nabla\cdot\bar h_n;
\end{equation}
For $n=N+1$, $\bar w_{N+1}$ denotes the finite energy solution of
\begin{equation}\nonumber
-\nabla\cdot a\nabla \bar w_{N+1}=\nabla\cdot\bar h_{N+1}+\bar f.
\end{equation}
By uniqueness of finite energy solutions of (\ref{4.5}), 
$\sum_{n=0}^{N+1}\nabla\bar w_n=\nabla \bar w$, so that
by the triangle inequality
\begin{equation}\label{6.5}
\left(\int_{B_1(x_0)}|\nabla\bar w|^2\right)^\frac{1}{2}\le
\sum_{n=0}^{N+1}\left(\int_{B_1(x_0)}|\nabla\bar w_n|^2\right)^\frac{1}{2}.
\end{equation}
In the sequel, we will estimate the contributions individually.

\medskip

We start with the intermediate $n=1,\ldots,N$: Since by construction of $\bar h_n$
and by (\ref{6.1}), $\bar w_n$ is $a$-harmonic in $B_{2^{n-1}}(x_0)$, we have by
Step 1 (the second inequality in (\ref{5.7}) and with the origin replaced by $x_0$), 
and the energy estimate for (\ref{6.1}):
\begin{align*}
\left(\int_{B_1(x_0)}|\nabla\bar w_n|^2\right)^\frac{1}{2}
&\lesssim
(\frac{1}{2^n})^\frac{d}{2}
\left(\int_{B_{2^{n-1}}(x_0)}|\nabla\bar w_n|^2\right)^\frac{1}{2}\\
&\lesssim
(\frac{1}{2^n})^\frac{d}{2}\left(\int_{B_{2^{n}}(x_0)}|\bar h|^2\right)^\frac{1}{2}.
\end{align*}
In the case of $n=0$, we obtain likewise by just the energy estimate
\begin{equation*}
\left(\int_{B_1(x_0)}|\nabla\bar w_0|^2\right)^\frac{1}{2}\lesssim\left(\int_{B_1(x_0)}|\bar h|^2\right)^\frac{1}{2}.
\end{equation*}
Hence in both cases we obtain thanks to (\ref{4.3}) that
\begin{equation}\label{6.3}
\left(\int_{B_1(x_0)}|\nabla\bar w_n|^2\right)^\frac{1}{2}\lesssim\frac{1}{R^{d+\alpha}}.
\end{equation}

\medskip

We finally turn to $n=N+1$ and obtain like for $n=1,\ldots,N$ by Step 1
\begin{equation*}
\left(\int_{B_1(x_0)}|\nabla\bar w_{N+1}|^2\right)^\frac{1}{2}\lesssim
(\frac{1}{R})^\frac{d}{2}
\left(\int_{B_{R}(x_0)}|\nabla\bar w_{N+1}|^2\right)^\frac{1}{2}.
\end{equation*}
But now we use the triangle inequality in form of
\begin{equation*}
\left(\int_{B_{R}(x_0)}|\nabla\bar w_{N+1}|^2\right)^\frac{1}{2}
\le
\left(\int_{B_{R}(x_0)}|\nabla\bar w|^2\right)^\frac{1}{2}
+\left(\int|\nabla(\bar w-\bar w_{N+1})|^2\right)^\frac{1}{2},
\end{equation*}
and the energy estimate for
$-\nabla\cdot a\nabla(\bar w-\bar w_{N+1})=\nabla\cdot(I(B_R(x_0))\bar h)$
to conclude
\begin{equation*}
\left(\int|\nabla (\bar w-w_{N+1})|^2\right)^\frac{1}{2}
\le\left(\int_{B_{R}(x_0)}|\bar h|^2\right)^\frac{1}{2}
\;\stackrel{(\ref{4.3})}{\lesssim}\;\frac{1}{R^{\frac{d}{2}+\alpha}}.
\end{equation*}
Collecting these estimates on the contribution of $n=N+1$, we obtain
\begin{equation}\label{6.4}
\left(\int_{B_1(x_0)}|\nabla\bar w_{N+1}|^2\right)^\frac{1}{2}
\lesssim
\frac{1}{R^{d+\alpha}}+
(\frac{1}{R})^\frac{d}{2}\left(\int_{B_{R}(x_0)}|\nabla\bar w|^2\right)^\frac{1}{2}.
\end{equation}
Inserting (\ref{6.3}) and (\ref{6.4}) into (\ref{6.5}), we obtain (\ref{6.6}).

\medskip

{\sc Argument for Step 7}. Note that by (\ref{4.6}) and (\ref{3.10}) we have in $B_1(x_0)$,
where the cut-off satisfies $\eta\equiv 1$,
\begin{equation}\nonumber
\nabla w=\nabla u-[\partial_iv(e_i+\nabla\phi_i)+(\phi_i-\fint_{B_1(x_0)}\phi_i)\nabla\partial_iv].
\end{equation}
In fact, the extra term $(\phi_i-\fint_{B_1(x_0)}\phi_i)\nabla\partial_iv$ is of
higher order:
\begin{equation}\nonumber
\left(\int_{B_1(x_0)}|(\phi_i-\fint_{B_1(x_0)}\phi_i)\nabla\partial_iv|^2\right)^\frac{1}{2}
\lesssim\frac{1}{R^{d+1}},
\end{equation}
which follows immediately from the assumption \eqref{T.1} and from the pointwise
estimate \eqref{2.1} on $v$.

\bigskip

{\sc Proof of Corollary \ref{C}}. Like for the proof of the Theorem \ref{T} we may assume $r_*=1$.
We consider $g$, $u$, and $v$ like in the statement of Theorem \ref{T} and note that we obtain
from (\ref{2.2}) the Green's function representation
\begin{equation}\nonumber
u(x)=-\int\nabla_yG(x,y)\cdot g(y)\dy,
\end{equation}
whereas (\ref{2.3}) yields
\begin{equation}\nonumber
v(x)=\int\nabla G_h(x-y)\cdot (e_i+\partial_i\phi(y)) g_i(y)\dy 
=\int\partial_jG_h(x-y)(e_j+\nabla\phi_j(y))\cdot g(y)\dy.
\end{equation}
By differentiation in $x$ this implies
\begin{align*}
\nabla u(x)&=-\int_{B_1(0)}\nabla_x\nabla_yG(x,y)g(y)\dy,\\
\partial_i v(x)&=\int_{B_1(0)}\partial_i\partial_jG_h(x-y)(e_j+\nabla\phi_j(y))\cdot g(y)\dy,
\end{align*}
so that (\ref{c1}) takes on the form
\begin{multline}
\Big(\int_{B_1(x_0)}\Big|\int_{B_1(0)}\Big(\nabla_x\nabla_yG(x,y)+\partial_i\partial_jG_h(x-y)\big(e_i+\nabla\phi_i(x)\big)\otimes\big(e_j+\nabla\phi_j(y)\big)\Big)g(y)\dy
                    \Big|^2\dx\Big)^\frac{1}{2}\\
\lesssim\frac{\ln|x_0|}{|x_0|^{d+\alpha}}\big(\int_{B_1(0)}|g|^2\dy\big)^\frac{1}{2}.\label{c2}
\end{multline}

\medskip

We now argue that in (\ref{c2}), we may replace $\partial_i\partial_jG_h(x-y)$ by $\partial_i\partial_jG_h(x)$.
Indeed, because of $|x_0|\ge 4$, $|x-x_0|\le 1$, and $|y|\le 1$, we have 
for the constant-coefficient Green's function 
$|\partial_i\partial_jG_h(x-y)-\partial_i\partial_jG_h(x)|\lesssim\frac{1}{|x_0|^{d+1}}$. In addition,
we have by the argument from Step 1 in the proof of Theorem \ref{T} and our assumption (\ref{T.1}) that
$\int_{B_1(x_0)}|e_i+\nabla\phi_i(x)|^2\dx\lesssim 1$ as well as $\int_{B_1(0)}|e_j+\nabla\phi_j(y)|^2\dy\lesssim 1$.
We therefore obtain by Cauchy-Schwarz' inequality in $y$
%
\begin{multline}
\Big(\int_{B_1(x_0)}\Big|\int_{B_1(0)}\big(\partial_i\partial_jG_h(x-y)-\partial_i\partial_jG_h(x)\big)\big(e_i+\nabla\phi_i(x)\big)\otimes\big(e_j+\nabla\phi_j(y)\big)g(y)\dy\Big|^2\dx\Big)^\frac{1}{2}\nonumber\\
\lesssim\frac{1}{|x_0|^{d+1}}\big(\int_{B_1(0)}|g|^2\dy\big)^\frac{1}{2}.
\end{multline}
Hence (\ref{c2}) upgrades to
\begin{multline}
\Big(\fint_{B_1(x_0)}\Big|\int_{B_1(0)}\Big(\nabla_x\nabla_yG(x,y)+\partial_i\partial_jG_h(x)\big(e_i+\nabla\phi_i(x)\big)\otimes\big(e_j+\nabla\phi_j(y)\big)\Big)g(y)\dy\Big|^2\dx\Big)^{\frac{1}{2}}\\
\lesssim\frac{\ln|x_0|}{|x_0|^{d+\alpha}}\big(\int_{B_1(0)}|g|^2\dy\big)^\frac{1}{2}.\label{c5}
\end{multline}
We now apply Lemma \ref{L} to the family $u(y)=\nabla_xG(x,y)-\partial_i\partial_jG_h(x)(e_i+\nabla\phi_i(x))
(y_j+\phi_j(y))$ of $\mathbb{R}^d$-valued maps defined for $y\in B_1(0)$
and parameterized by the point $x\in B_1(x_0)$, i.e., we interchanged the role of $x$ and $y$. We note that these maps are component-wise $a$-harmonic
on $\{|y|\le 1\}$ because of (\ref{i2}) and because of the $y$-derivative of (\ref{c4}) in conjunction with the symmetry of the Green's function, which follows from symmetry of $A$ (for nonsymmetric $A$ one would apply Lemma \ref{L} for the adjoint problem). The ensemble average on this family is given by the spatial average $\langle\cdot\rangle=\fint_{B_1(x_0)}\cdot \dx$. 
The role of the linear functionals in the statement of Lemma \ref{L} is played by 
$Fu:=\int_{B_1(0)}\nabla_y u(y) g(y)\dy$, where we restrict to $g$ that are normalized
$(\int_{B_1(0)}|g|^2\dy)^\frac{1}{2}=1$. Hence we learn from Lemma \ref{L} that (\ref{c5}) implies
\begin{align*}
\Big(\!\fint_{B_1(x_0)}\int_{B_\frac{1}{2}(0)}\!\big|\nabla_x\nabla_yG(x,y)+\partial_i\partial_jG_h(x)(e_i+\nabla\phi_i(x))\otimes(e_j+\nabla\phi_j(y))\big|^2\!\dy\dx\Big)^\frac{1}{2}\lesssim\frac{\ln|x_0|}{|x_0|^{d+\alpha}}.
\end{align*}
By the preceding argument, we may substitute $\partial_i\partial_jG_h(x)$ 
again by $\partial_i\partial_jG_h(x-y)$.

\bigskip

{\sc Proof of Lemma \ref{L}}.
Using translation and scaling invariance, an elementary covering argument shows that is enough to establish
(\ref{b10}) with radius $\frac{R}{2}$ replaced by $\frac{R}{2\sqrt{d}}$. Therefore,
it suffices to show the result with
the inner ball $\{|x|<\frac{R}{2\sqrt{d}}\}$ replaced by the cube $(-\frac{R}{2},\frac{R}{2})^d$ and the outer
ball $\{|x|<R\}$ in (\ref{b9}) by the cube $(-R,R)^d$. By scale invariance, we may reduce to
$(-\frac{\pi}{4},\frac{\pi}{4})^d$ and $(-\frac{\pi}{2},\frac{\pi}{2})^d$, respectively.
We thus will show that
\begin{equation}\label{b5}
\langle\int_{(-\frac{\pi}{4},\frac{\pi}{4})^d}|\nabla u|^2\dx\rangle\lesssim
\sup_{F}{\langle |F u|^2\rangle},
\end{equation}
where the supremum is taken over all linear functionals $F$ satisfying
\begin{equation}\label{b3}
|Fu|^2\le\int_{(-\frac{\pi}{2},\frac{\pi}{2})^d}|\nabla u|^2\dx.
\end{equation}

\medskip

The proof essentially amounts to a generalization of Caccioppoli's estimate
for an $a$-harmonic function $u$. Recall that the latter states
\begin{equation}\label{b1}
\int_{(-\frac{\pi}{4},\frac{\pi}{4})^d}|\nabla u|^2\dx
\lesssim\inf_{c\in\mathbb{R}}\int_{(-\frac{\pi}{2},\frac{\pi}{2})^d}|u-c|^2\dx,
\end{equation}
which we may re-express in terms of the Fourier cosine series
\begin{equation}\label{b4}
{\mathcal F}u(k):={\textstyle\sqrt\frac{2}{\pi^d}}
\int_{(-\frac{\pi}{2},\frac{\pi}{2})^d} u(x)\Pi_{i=1}^d\cos(k_ix_i)\dx
\quad\mbox{for}\;k\in\mathbb{Z}^d\setminus\{0\}
\end{equation}
as
$\int_{(-\frac{\pi}{4},\frac{\pi}{4})^d}|\nabla u|^2\dx
\lesssim\sum_{k\in\mathbb{Z}^d\setminus\{0\}}|{\mathcal F}u(k)|^2$.
The generalization of (\ref{b1}) we need is that for any even $n\in\mathbb{N}$ we have
\begin{equation}\label{b2}
\int_{(-\frac{\pi}{4},\frac{\pi}{4})^d}|\nabla u|^2\dx
\lesssim\sum_{k\in\mathbb{Z}^d\setminus\{0\}}\frac{1}{|k|^{2n}}|{\mathcal F}u(k)|^2,
\end{equation}
which amounts to replace the $L^2$-norm on the r.\ h.\ s.\ of (\ref{b1}) by
the negative $\dot H^{-n}$-norm. For the remainder of the proof, $\lesssim$
will mean up to a constant also depending on $n$; but this will not
matter, since we will presently fix an $n$ in terms of $d$.

\medskip

Before we give the argument for (\ref{b2}), let us argue how to conclude.
We first note that for $k\in\mathbb{Z}^d\setminus\{0\}$ the linear functional
\begin{equation}\nonumber
F_k u:=|k|{\mathcal F}u(k)
\end{equation}
has the boundedness property (\ref{b3}).
Indeed, we obtain from integration by parts in (\ref{b4})
\begin{equation}\nonumber
{\mathcal F}u(k)
=-\frac{1}{k_1}{\textstyle\sqrt\frac{2}{\pi^d}}
\int_{(-\frac{\pi}{2},\frac{\pi}{2})^d}\partial_1u(x)\sin(k_1x_1)\Pi_{i=2}^d\cos(k_ix_i)\dx,
\end{equation}
so that $k_1^2|{\mathcal F}u(k)|^2\le\int_{(-\frac{\pi}{2},\frac{\pi}{2})^d}(\partial_1u)^2(x)\dx$.
Hence, after taking the ensemble average, we may reformulate (\ref{b2}) as
\begin{equation}\nonumber
\langle\int_{(-\frac{\pi}{4},\frac{\pi}{4})^d}|\nabla u|^2\dx\rangle
\lesssim\sum_{k\in\mathbb{Z}^d\setminus\{0\}}\frac{1}{|k|^{2(n+1)}}\langle|F_ku|^2\rangle.
\end{equation}
Now picking $n\in\mathbb{N}$ with $n>\frac{d}{2}-1$ so that
$\sum_{k\in\mathbb{Z}^d\setminus\{0\}}\frac{1}{|k|^{2(n+1)}}\lesssim 1$, we obtain (\ref{b5}).

\medskip

We now turn to the argument for (\ref{b2}) and introduce the abbreviation $\|\cdot\|$
for the $L^2((-\frac{\pi}{2},\frac{\pi}{2})^d)$-norm. The main ingredient is the following interpolation
inequality for any function $v$ of zero spatial average
\begin{equation}\label{b7}
\|\eta^n v\|
\lesssim\|\eta^{n+1}\nabla v\|^\frac{n}{n+1}
\Big(\!\!\!\sum_{k\in\mathbb{Z}^d\setminus\{0\}}\frac{1}{|k|^{2n}}|{\mathcal F}v(k)|^2\Big)^\frac{1}{2(n+1)}
+\Big(\sum_{k\in\mathbb{Z}^d\setminus\{0\}}\frac{1}{|k|^{2n}}|{\mathcal F}v(k)|^2\Big)^\frac{1}{2},
\end{equation}
where $\eta$ is a cut-off function for $(-\frac{\pi}{4},\frac{\pi}{4})^d$ in
$(-\frac{\pi}{2},\frac{\pi}{2})^d$ with
\begin{equation}\label{b8}
|\nabla\eta|\lesssim1.
\end{equation}
Note that (\ref{b7}) couples the degree of negativity of the r.\ h.\ s.\ norm to
the degree of degeneracy of the cut-off $\eta^n$.
If we plug the standard Caccioppoli estimate in its refined form of
\begin{equation}\nonumber
\|\eta^{n+1}\nabla u\|
\lesssim\inf_{c\in\mathbb{R}}\|(u-c)\nabla\eta^{n+1}\|
\stackrel{(\ref{b8})}{\lesssim}\inf_{c\in\mathbb{R}}\|\eta^{n}(u-c)\|
\end{equation}
into (\ref{b7}) for $v=u-c$ and use Young's inequality, we obtain (\ref{b2}).

\medskip
In preparation of its proof, we rewrite (\ref{b7}) without Fourier transform,
appealing to the representation of the Laplacian $-\triangle_N$ with Neumann boundary
conditions through the Fourier cosine series by ${\mathcal F}(-\triangle_N)w(k)=|k|^2{\mathcal F}w(k)$:
\begin{equation*}
\Big(\sum_{k\in\mathbb{Z}^d\setminus\{0\}}\frac{1}{|k|^{2n}}|{\mathcal F}v(k)|^2\Big)^\frac{1}{2}
=\|w\|\quad\mbox{where}\quad(-\triangle_N)^\frac{n}{2}w=v.
\end{equation*}
For (\ref{b7}) it thus suffices to show for arbitrary function $w$
\begin{equation*}
\|\eta^{n}\triangle^\frac{n}{2}w\|
\lesssim\|\eta^{n+1}\nabla\triangle^\frac{n}{2}w\|^\frac{n}{n+1}
\|w\|^\frac{1}{n+1}+\|w\|.
\end{equation*}
By iterated application of Young's inequality, it is easily seen that this family of interpolation
estimates indexed by even $n$ follows from the following two-tier family of
interpolation inequalities index by $m\in\mathbb{N}$
\begin{equation*}
\|\eta^{2m}\triangle^m w\|
\lesssim\|\eta^{2m+1}\nabla\triangle^mw\|^\frac{1}{2}
\|\eta^{2m-1}\nabla\triangle^{m-1} w\|^\frac{1}{2}
+\|\eta^{2m-1}\nabla\triangle^{m-1} w\|
\end{equation*}
and
\begin{equation*}
\|\eta^{2m-1}\nabla \triangle^{m-1} w\|
\lesssim\|\eta^{2m}\triangle^{m}w\|^\frac{1}{2}
\|\eta^{2m-2}\triangle^{m-1} w\|^\frac{1}{2}
+\|\eta^{2m-2}\triangle^{m-1} w\|.
\end{equation*}
Obviously, this two-tier family reduces to the two estimates
\begin{align*}
\|\eta^{2m}\triangle v\|
&\lesssim\|\eta^{2m+1}\nabla\triangle v\|^\frac{1}{2}
\|\eta^{2m-1}\nabla v\|^\frac{1}{2}
+\|\eta^{2m-1}\nabla v\|,\\
\|\eta^{2m-1}\nabla v\|
&\lesssim\|\eta^{2m}\triangle v\|^\frac{1}{2}
\|\eta^{2m-2} v\|^\frac{1}{2}
+\|\eta^{2m-2} v\|,
\end{align*}
which by Young's inequality follow from
\begin{align*}
\|\eta^{2m}\triangle v\|
&\lesssim(\|\eta^{2m+1}\nabla\triangle v\|+\|\eta^{2m}\triangle v\|)^\frac{1}{2}
\|\eta^{2m-1}\nabla v\|^\frac{1}{2},\\
\|\eta^{2m-1}\nabla v\|
&\lesssim(\|\eta^{2m}\triangle v\|+\|\eta^{2m-1}\nabla v\|)^\frac{1}{2}
\|\eta^{2m-2} v\|^\frac{1}{2}.
\end{align*}
Thanks to (\ref{b8}), these two last estimates immediately follow
from integration by parts (the cut-off $\eta$ suppresses boundary terms),
the Cauchy-Schwarz and the triangle inequalities.

\providecommand{\bysame}{\leavevmode\hbox to3em{\hrulefill}\thinspace}
\providecommand{\MR}{\relax\ifhmode\unskip\space\fi MR }
\providecommand{\MRhref}[2]{%
  \href{http://www.ams.org/mathscinet-getitem?mr=#1}{#2}
}
\providecommand{\href}[2]{#2}

\end{document}